\documentclass[a4paper,11pt]{article}
\usepackage{}
\usepackage{mathrsfs}
\usepackage{amssymb}
\usepackage{amsmath}
\usepackage{amsfonts}
\usepackage{amsmath}
\usepackage{graphicx}
\numberwithin{equation}{section}

\newtheorem{lemma2}{Lemma}[section]

\newtheorem {remark2}{Remark}[section]

\newtheorem {theorem3}{Theorem}[section]
\newtheorem {remark3}{Remark}[section]
\newtheorem {lemma3}{Lemma}[section]

\newtheorem {theorem4}{Theorem}[section]
\newtheorem {corollary4}{Corollary}[section]
\newtheorem {remark4}{Remark}[section]

\newtheorem {theorem5}{Theorem}[section]
\newtheorem {lemma5}{Lemma}[section]

\newtheorem {remark5}{Remark}[section]
\begin{document}

\title{\textbf{The Well-posedness and Blow-up rate of Solution for the Generalized
Zakharov equations with Magnetic field in $\mathbb{R}^d$}}
\author{
Xinglong Wu \footnote{Email: wxl8758669@aliyun.com}\\
Wuhan Institute of Physics and Mathematics,\\ Chinese Academy of Sciences,
 Wuhan 430071, P. R. China \smallskip
\\
Boling Guo\\
Institute of Applied Physics and Computational Mathematics\\
 Beijing, 100088, P. R. China %\footnote{E-mail: gbl@mail.iapcm.ac.cn}
 }
\date{}
\maketitle

\begin{abstract}
The present paper is devoted to the study of the well-posedness and
the lower bound of blow-up rate to the Cauchy problem of the
generalized Zakharov(GZ) equations with magnetic field in $\mathbb{R}^d$,
$d\geq1$. The work of well-posedness of the GZ system  bases on the
local well-posedness theory in \cite{GTV}. At first,
the existence, uniqueness and continuity of solution to the
GZ system with magnetic field in
$\mathbb{R}^d$ is proved. Next, we establish the lower bound of blow-up rate
of blow-up solution in sobolev spaces to the GZ system,
which is almost a critical index. Finally, we obtain the long time behavior of global solution, whose $H^k$-norm  grows at $k$-exponentially in time.


\end{abstract}
\section{Introduction}
In this paper, we consider the Cauchy problem of the generalized
Zakharov system with magnetic field in $\mathbb{R}^d$ as follows
\begin{equation}
\left\{\begin{array}{ll}i\partial_{t}E+\Delta E -nE+i(E\otimes
B)=0,&(t,x)\in\mathbb{R}^+\times\mathbb{R}^d,
\\\frac{1}{c_0^2}\partial_{tt}n-\Delta (n+|E|^2)=0,\\
\Delta B-i\eta\nabla\times\nabla\times(E\otimes \overline{E})+A=0,
\\E(0,x)=E_0(x),&x\in\mathbb{R}^d,\\
(n(0,x),\partial_tn(0,x)) =(n_{0}(x),n_{1}(x)),\end{array}\right.
\end{equation}
where $c_0>0$ is a constant, $E(t,x)$ denotes a vector valued
function from $\mathbb{R}^+\times\mathbb{R}^d$ into $\mathbb{C}^d$.
$n(t,x)$ is a function from $\mathbb{R}^+\times \mathbb{R}^d$ to
$\mathbb{R}$, $B(t,x)$ is a vector valued function from
$\mathbb{R}^+\times \mathbb{R}^d$ into $\mathbb{C}^d$, and $A$ has
the following two
form: \\
$(A_1)\qquad A=\beta B,\quad \beta$ is a nonpositive constant;\\
$(A_2)\qquad A=-\gamma\frac{\partial}{\partial
t}\int_{\mathbb{R}^d}\frac{B(t,y)}{|x-y|^2}dy.$\\
 The system (1.1) describes the spontaneous generation of a magnetic
 field in a cold plasma (case $A_1$) or in a hot plasma (case $A_2$)
 \cite{KSH}. $E$ denotes the slowly varying complex amplitude of the
high-frequency electric field, $n(t, x)$ represents the
 fluctuation of the electron density from its equilibrium, $B$ is
 the self-generated magnetic field. $i^2 =-1$, constant $\eta>0$,
 $\overline{E}$ denotes the complex conjugate of $E$, and $\otimes$ means the exterior
product of vector-valued functions.
\par
If we neglect the magnetic field $B$, system (1.1) becomes the
classical Zakharov equation
\begin{equation}
\left\{\begin{array}{ll}i\partial_{t}E+\Delta E
-nE=0,(t,x)\in\mathbb{R}^+\times\mathbb{R}^d,
\\\frac{1}{c_0^2}\partial_{tt}n-\Delta (n+|E|^2)=0,
\\E(0,x)=E_0(x),x\in\mathbb{R}^d,
\\(n(0,x),\partial_tn(0,x)) =(n_{0}(x),n_{1}(x)),\end{array}\right.
\end{equation}
which describes the propagation of Langmuir wave \cite{Za}. The
Cauchy problem of Eq.(1.2) was established by several authors
\cite{Ad, Ad1, B, GTV, KPV, OT, Su}. Such as, the local
well-posedness was obtain in spaces $H^k\times H^l\times
H^{l-1}$\cite{GTV} for any dimensions $d$, C. Sulem and P.L. Sulem
\cite{Su} proved the global existence of a weak solution in two and
three dimensions for the small initial data. With the same
assumptions, they also got the existence and uniqueness of the
smooth solution $(E,n)\in \mathcal {C}([0,T[;H^m)\times\mathcal
{C}([0,T[;H^{m-1}),m\geq3$. Moreover, the solution was global in one
dimension, and the global solution can be extended in two dimensions
with small initial data \cite{Ad}. Numerical simulations strongly
suggest a finite blow-up time for some initial data, and global
solution of the small initial data  can be numerically verified by
Papanicolaou, C. Sulem, P. L. Sulem, Wang, and Landman \cite{LPSS,
PSSW}. By constructing a family of blow-up solutions of the
following form
\begin{equation*}
\left\{\begin{array}{ll}E(t,x)=\frac{\omega}{T-t}e^{i(\theta+\frac{|x|^2-4\omega^2}{4(-T+t)})}
P\left(\frac{x\omega}{T-t}\right),
\\n(t,x)=\left(\frac{\omega}{T-t}\right)^2N\left(\frac{x\omega}{(T-t)}\right),\end{array}\right.
\end{equation*}
where $\omega>\omega_0,\theta\in\mathbb{R},$ and
$$P(x)=P(|x|),N(x)=N(|x|), \Delta P-P=NP,$$
$$\frac{1}{(c_0\omega)^2}(r^2N_{rr}+6rN_r+6N)-\Delta N=\Delta P^2,$$
with $r=|x|, \Delta w=w_{rr}+\frac{1}{r}w_r$, L. Glangetas and F.
Merle \cite{GM} proved the existence of self-similar blow-up
solutions to the Hamiltonian case of Eq.(1.2) in two dimensions.
i.e.
\begin{equation}
\left\{\begin{array}{ll}i\partial_{t}E+\Delta E -nE=0,
\\
\partial_tn+div v=0,\quad(t,x)\in\mathbb{R}^+\times\mathbb{R}^2,
\\\frac{1}{c_0^2}\partial_{t}v+\nabla
(n+|E|^2)=0,
\\E(0,x)=E_0(x),\quad x\in\mathbb{R}^2,
\\(n(0,x),\partial_tn(0,x)) =(n_{0}(x),n_{1}(x)),\end{array}\right.
\end{equation}
more results of Eq.(1.3) can be found in \cite{GM1, Me, Me1}. In
fact, the existence and uniqueness of global solution is open
problem in $d\geq3$. It is interested to recall the situation in the
case $c_0=\infty$, that is the Zakharov equations reduce to the cube
nonlinear Schr\"{o}dinger equation \cite{CW, G, LPSS, Na}
$$i\partial_tu+\Delta u=-|u|^2u.$$
\par
Returning to the generalized Zakharov system (1.1) with magnetic
field. We consider the system (1.1) in the Hamiltonian case, i.e.
 $$\partial_tn(t,x)=-\Delta w(t,x)=-div V(t,x),$$ then system (1.1) can
be written in the form \cite{La}
\begin{equation}
\left\{\begin{array}{ll}i\partial_{t}E+\Delta E-nE+i(E\otimes B)=0,
\\
\partial_tn+div V=0,\quad(t,x)\in\mathbb{R}^+\times\mathbb{R}^d,
\\\frac{1}{c_0^2}\partial_{t}V+\nabla
(n+|E|^2)=0,\\
\Delta B-i\eta\nabla\times\nabla\times(E\otimes \overline{E})+A=0,
\\E(0,x)=E_0(x),\quad x\in\mathbb{R}^d,
\\(n(0,x),V(0,x)) =(n_{0}(x),V_0(x)).\end{array}\right.
\end{equation}
In 1995, by the conservation laws of Eq.(1.4) in the case $(A_1)$
\begin{equation}H_1(E):=\|E\|_{L^2}=\|E_0\|_{L^2}=H_1(E_0),\end{equation}
\begin{equation}\begin{split}H_2(E,n,V,B)&:=\|\nabla
E\|_{L^2}^2+\frac{1}{2}\|n\|_{L^2}^2+\frac{1}{2c_0^2}\|V\|_{L^2}^2+\int_{\mathbb{R}^d}n|E|^2dx\\&-\frac{i}{2}
\int_{\mathbb{R}^d}
B(E\otimes\overline{E})dx=H_2(E_0,n_0,V_0,B_0),\end{split}\end{equation}
C. Laurey \cite{La} got the global existence of weak solution
$(E,n,B)\in L^\infty(\mathbb{R}^+;H^1)\times
L^\infty(\mathbb{R}^+;L^2)\times L^\infty(\mathbb{R}^+;L^2)$ to
Eq.(1.4) in the case $(A_1)$ with the small initial data
$(E_0,n_0,B_0)\in H^1\times L^2\times L^2$. As the initial data
$(E_0,n_0,B_0)\in H^{s+1}\times H^s\times H^{s+1},s>\frac{d}{2},
d=2,3$, he established the local existence and uniqueness of a
strong solution $(E,n,B)\in L^\infty([0,T[;H^{s+1})\times
L^\infty([0,T[;H^s)\times L^\infty([0,T[;H^{s+1})$ to system (1.1)
in the case $(A_1)$ and $(A_2)$, for some $T>0$. If $d=2,s=2$, in
the case $(A_1)$, the smooth solution was global with the small
initial data. Recently, similar to \cite{GM}, in two dimensions,
Gan, Guo and Huang \cite{GGH} constructed a family of blow-up
solutions and proved the existence of self-similar blow-up solution
to Eq.(1.4) in the case $(A_1)$.
 \par
A natural problem of system (1.1) is to establish the global
solution or construct the blow-up solution in dimensions $d\geq3$.
In this paper, at first, for the  generalized  Zakharov equation
(1.1) with magnetic field in $\mathbb{R}^d$,  similar to \cite{GTV}
the well-posedness is obtained in spaces
$$(E,B,n,\partial_tn)\in\mathcal
{C}([0,T[;H^k)\times\mathcal {C}([0,T[;H^k)\times\mathcal
{C}([0,T[;H^l)\times\mathcal {C}([0,T[;H^{l-1}),$$ if the initial
data $(E_0,B_0,n_0,n_1)$ belongs to $H^k\times H^k\times H^l\times
H^{l-1}$. The difficult is how to deal with the nonlinearity in
system $(1.1)_1$ and system $(1.1)_3$ with $(A_1)$ or $(A_2)$ is
decisive. If the solution blows up in finite time $T^*$, then the
lower bound for the blow-up rate of the blow-up solution to system
(1.1) satisfies
$$\|E(t)\|_{H^l}+\|n(t)\|_{H^l}+\|\partial_tn\|_{H^{l-1}}>C\frac{1}{(T^*-t)^{\frac{3}{4}\vartheta_l-\epsilon}},$$
where $\vartheta_l=\frac{1}{4}(2l+4-d)$ and $\epsilon$ is a any
positive constant. Moreover, for equation (1.2) in 3D, we have the
following lower bound
$$\|E(t)\|_{H^l}+\|n(t)\|_{H^l}+\|\partial_tn\|_{H^{l-1}}>C\frac{1}{(T^*-t)^{\vartheta_l-\epsilon}},$$which
almost up to the bound of the following asymptotic self-similar
blow-up solution to equation (1.2)
\begin{equation*}
\left\{\begin{array}{ll}E(t,x)=\frac{1}{T-t}
P\left(\frac{|x|}{(T-t)^{1/2}}\right)+i\frac{1}{T-t}Q\left(\frac{|x|}{(T-t)^{1/2}}\right),
\\n(t,x)=\frac{1}{T-t}N\left(\frac{|x|}{(T-t)^{1/2}}\right),\end{array}\right.
\end{equation*}
where $P(x)=P(|x|),Q(x)=Q(|x|),N(x)=N(|x|)$, and $(P,N)$ satisfies
the ODEs
\begin{equation*} \left\{\begin{array}{ll}\Delta P+\frac{1}{2}rQ_r+Q=NP,\\
\Delta Q+\frac{1}{2}rP_r+P=NQ,\\
\frac{1}{4}r^2N_{rr}+\frac{7}{4}rN_r+2N)=\Delta
(P^2)+\Delta(Q^2).\end{array}\right.
\end{equation*}
\par
 The remainder of this paper is organized as follows. In Section 2,
we first recall the definition of the weighted Bourgain spaces, some
important lemmas, and present the proof's frame for the
well-posedness of the Cauchy problem to system (1.1). In Section 3,
the local well-posedness of the Cauchy problem to system (1.1) with
the magnetic field $B$ satisfying the case $(A_1)$, $(A_2)$ is
established in $H^k\times H^l\times H^{l-1}$. Next, in Section 4, we
derive the lower bound of blow-up rate of blow-up solution in
Sobolev spaces to system (1.1) and Eq.(1.2), which is almost a critical index.  In Section 5, we obtain the global solution to system (1.1) with the small initial data, and the $H^k$-norm of solution grows at $k$-exponentially in time.\\
%Notation: Let $\widetilde{L}^\rho_T(\dot{B}_{p,r}^\sigma)$ denote
%the set of functions $u$ such that
%$$\|u\|_{\widetilde{L}^\rho_T(\dot{B}_{p,r}^\sigma)}:=\left\|\left(2^{k\sigma}\|\dot{\Delta}_ku\|_{L_T^\rho(L^p)}\right)_{k\in\mathbb{Z}}\right\|_{l^r},$$
% for $\sigma\in\mathbb{R},T>0$, and $(p,r,\rho)\in[1,\infty]^3$. Then
% via the Minkowski inequality, we have
% $$\|u\|_{\widetilde{L}^\rho_T(\dot{B}_{p,r}^\sigma)}\leq C\|u\|_{L^\rho_T(\dot{B}_{p,r}^\sigma)}\qquad \text{if}\quad 1\leq\rho\leq r,$$
%otherwise,
%$$\|u\|_{L^\rho_T(\dot{B}_{p,r}^\sigma)}\leq C\|u\|_{\widetilde{L}^\rho_T(\dot{B}_{p,r}^\sigma)}\qquad \text{if}\quad 1\leq r\leq\rho $$
%If for any $T>0$, $u\in\widetilde{L}^\rho_T(\dot{B}_{p,r}^\sigma)$,
%then we have $u\in\widetilde{L}^\rho(\dot{B}_{p,r}^\sigma)$, where
%$$\|u\|_{\widetilde{L}^\rho(\dot{B}_{p,r}^\sigma)}:=\left(\sum_{k\in\mathbb{Z}}(2^{k\sigma}\|\dot{\Delta}_ku
%\|_{L^\rho(\mathbb{R}^+;L^p)})^r\right)^{1/r}.$$
\section{The Preliminary}
In this subsection, for the convenience of the readers, we recall
the Bourgain method, in order to make this paper self-contained and
to locate exactly the required nonlinear estimated, which are come
from the paper\cite{GTV}, only make some little modification for our
target. These important lemmas will be used repeatedly throughout
this paper.
\par In order to deal with the second wave equation of
system (1.1), without loss of generality, let $c_0=1$, we split $n$
into its positive and negative frequency parts as
\begin{equation}
\varphi_\pm=n\pm i\Lambda^{-1}\partial_tn,
\end{equation}
where the operator $\Lambda=(1-\Delta)^{\frac{1}{2}}$. Then one can
easily check that
$$(i\partial_t\mp \Lambda)\varphi_\pm=\mp\Lambda^{-1}(\partial_{tt}-\Delta)n\mp\Lambda^{-1}n.$$
Therefore, system $(1.1)_2$ is equivalent to
\begin{equation}
i\partial_t\varphi_\pm=\pm\Lambda\varphi_\pm\mp\Lambda^{-1}(\Delta|E|^2)\mp\frac{1}{2}\Lambda^{-1}(\varphi_++\varphi_-).
\end{equation}
The generalized Zakharov system (1.1) then takes the following form
\begin{equation}
\left\{\begin{array}{ll}i\partial_{t}E+\Delta E
-\frac{1}{2}(\varphi_++\varphi_-)E+i(E\otimes B)=0,
\\i\partial_t\varphi_\pm=\pm\Lambda\varphi_\pm\mp\Lambda^{-1}(\Delta|E|^2)\mp\frac{1}{2}\Lambda^{-1}(\varphi_++\varphi_-),\\
\Delta B-i\eta\nabla\times\nabla\times(E\otimes \overline{E})+A=0,
\end{array}\right.
\end{equation}
with the initial data $(E_0,\varphi_{0\pm},B_0)=(E_0,n_0\pm
i\Lambda^{-1}n_1,B_0)$.
\par
Define the semigroup $S(t)=e^{it\Delta}$,
$W^+(t)=e^{-it(1-\Delta)^{\frac{1}{2}}}$ and $W^-(t)=e^{-it(1-\Delta)^{\frac{1}{2}}}$,  the solution
$(E,\varphi_\pm)$ of the Cauchy problem to Eq.(2.3) is rewritten in
a standard way as integral equation
\begin{equation}
E=S(t)E_0-i\int_0^tS(t-\tau)[\frac{1}{2}(\varphi_++\varphi_-)E-i(E\otimes
B)](\tau)d\tau
\end{equation}
and
\begin{equation}
\varphi_\pm=W^{\pm}(t)\varphi_{0\pm}+i\int_0^tW^{\pm}(t-\tau)[\mp\Lambda^{-1}
(\Delta|E|^2)\mp\frac{1}{2}\Lambda^{-1}(\varphi_++\varphi_-)](\tau)d\tau.
\end{equation}
In order to use function space norms defined in terms of the space
time Fourier transform of solutions $(E,\varphi_\pm)$ in the context
on finite time interval $[-T,T]$, we introduce an even time cut-off
function $\psi\in\mathcal {C}_0^\infty$ satisfying
\begin{equation*}
\left\{\begin{array}{ll}\psi(t)=1,& if\quad|t|\leq1,
\\0\leq\psi(t)\leq1,& if\quad|t|\geq0,\\
\psi(t)=0,& if\quad|t|\geq2.
\end{array}\right.
\end{equation*}
Denote $\psi_{T}(t)=\psi(\frac{t}{T})$. Consider the cut-off
equation
\begin{equation}
E=\psi_1(t)S(t)E_0-i\psi_T(t)\int_0^tS(t-\tau)f_1(\tau)d\tau
\end{equation}
and
\begin{equation}
\varphi_\pm=\psi_1(t)W^\pm(t)\varphi_{0\pm}+i\psi_T(t)\int_0^tW^\pm(t-\tau)\left[f(\tau)\mp\frac{1}{2}\Lambda^{-1}\psi_{2T}
(\varphi_++\varphi_-)(\tau)\right]d\tau,
\end{equation}
where
$f_1=\psi_{2T}^2[\frac{1}{2}(\varphi_++\varphi_-)E]-i\psi_{3T}^3[(E\otimes
B)]$, and $f=\mp\Lambda^{-1} \psi_{2T}^2(\Delta|E|^2)$. One can
easily check that (2.6), (2.7) is actually identical with (2.4),
(2.5) respectively on Spp$\psi_T$, if $|t|\leq T\leq1$. For convenience, let $W$ denotes $W^+$ and $W^-$ in this paper.
\par
Introducing the space-time
weighted Bourgain spaces with norms respectively given by
$$\|u\|_{X_S^{s,b}}:=\|<\xi>^s<\tau+|\xi|^2>^b\hat{u}(\tau,\xi)\|_{L^2_{\tau,\xi}},$$
$$\|v\|_{X_W^{s,b}}:=\|<\xi>^s<\tau+|\xi|>^b\hat{v}(\tau,\xi)\|_{L^2_{\tau,\xi}},$$
where $<\cdot>=<1+|\cdot|^2>^{1/2}$. In a similar way, define the
$Y_S^k$ and $Y_W^k$ space with the norm
$$\|u\|_{Y_S^k}:=\|<\xi>^k<\tau+|\xi|^2>^{-1}\hat{u}(\tau,\xi)\|_{L^2_{\xi}(L_\tau^1)},$$
and$$\|u\|_{Y_W^k}:=\|<\xi>^k<\tau+|\xi|>^{-1}\hat{u}(\tau,\xi)\|_{L^2_{\xi}(L_\tau^1)}.$$
In order to solve the Cauchy problem of the generalized Zakharov
system (1.1) in the form of the integral equation (2.6), (2.7) by
the contraction mapping theorem in the space $X_{S}^{s,b},
X_{W}^{k,l}$, similar to the method in \cite{GTV}, we recall some
important lemmas which make some little modification for our target.
\begin{lemma2} Assume $s,b\in\mathbb{R}$. Then we have
$$\|\psi_1S(t)u_0\|_{X_{S}^{s,b}}\leq\|\psi_1\|_{H_t^b}\|u_0\|_{H_x^s},$$and
$$\|\psi_1W(t)u_0\|_{X_{W}^{s,b}}\leq\|\psi_1\|_{H_t^b}\|u_0\|_{H_x^s}.$$
Moreover, define $\psi(S\ast_
Rf)=\psi\int_0^tS(t-\tau)f(\tau)d\tau$, if $-\frac{1}{2}<b'\leq0\leq
b\leq b'+1$, and $T\leq1$, then we have
\begin{equation}\|\psi_T(S\ast_
Rf)\|_{X_{S}^{s,b}}\leq CT^{1-b+b'}\|f\|_{X_{S}^{s,b'}},
\end{equation}
and
\begin{equation}\|\psi_T(W\ast_
Rf)\|_{X_{W}^{s,b}}\leq CT^{1-b+b'}\|f\|_{X_{W}^{s,b'}}.
\end{equation}
\end{lemma2}
\begin{remark2} While we deal with the (2.6) and (2.7), the
following inequality will be used in order to get the positive power
of $T$,$$\|\psi_Tu\|_{X_{S}^{s,b}}\leq
CT^{-b+\frac{1}{q}}\|u\|_{X_{S}^{s,b}},
$$
where $s\in\mathbb{R}$, $b\geq0,\; 2\leq q$ and $1<bq$. if
$X_W^{s,b}$ take place of $X_S^{s,b}$, the above result is also
right.
\end{remark2}

\begin{lemma2} Let $b_0>\frac{1}{2}$ and $0\leq\gamma\leq1$. Assume
$ a,a_1,a_2\geq0$ satisfy
$$(1-\gamma)a<b_0,$$
\begin{equation}(1-\gamma)\max(a,a_1,a_2)\leq b_0\leq(1-\gamma)(a+a_1+a_2),
\end{equation}
\begin{equation}m\geq\frac{d}{2}+1-(1-\gamma)(a+a_1+a_2)/b_0\geq0\end{equation}
with strict inequality in (2.11L) if equality holds in (2.10R) or if
$a_1=0$. Let $v,v_1,v_2\in L^2$ such that $\mathcal
{F}^{-1}(\langle\tau+|\xi|\rangle^{\gamma
a}\hat{v}),\mathcal{F}^{-1}(\langle\tau+|\xi_i|\rangle^{\gamma
a_i}\hat{v_i}),i=1,2$ have support in $|t|\leq CT$. Then we deduce
the following estimates
\begin{equation}
\int\frac{|\hat{v}\hat{v}_1\hat{v}_2|}{\langle\sigma\rangle^a\langle\sigma_1\rangle^{a_1}\langle\sigma_2\rangle^{a_2}\langle\xi\rangle^m}\leq
CT^{\gamma(a+a_1+a_2)}\|v\|_{L^2}\|v_1\|_{L^2}\|v_2\|_{L^2},
\end{equation}
\begin{equation}
\int\frac{\hat{v}\hat{v}_1\hat{v}_2}{\langle\sigma\rangle^a\langle\sigma_1\rangle^{a_1}\langle\sigma_2\rangle^{a_2}\langle\xi_2\rangle^m}
\leq CT^{\gamma(a+a_1+a_2)}\|v\|_{L^2}\|v_1\|_{L^2}\|v_2\|_{L^2},
\end{equation}
where $\sigma=\tau+|\xi|,\quad\sigma_i=\tau+|\xi_i|^2,i=1,2$.
\end{lemma2}
By virtue of Lemma 2.2, we obtain Lemma 3.4 in \cite{GTV} with the
positive power $\theta_1$ of $T$, the $\theta_1$ can be up to
$\gamma(b+b_1+c_1)$, we also get Lemma 3.5, Lemma 3.6, Lemma 3.7 of
\cite{GTV}, the positive power of $T$ is $\theta_2=\gamma(c+2b_1)$,
$\theta_3=\gamma(b+1+b_1)$, $\theta_4=\gamma(1+2b_1)$ respectively.
\par
Similar to the proof of Proposition 3.1 in \cite{GTV}, taking
advantage of Lemma 2.1 and Lemma 2.2 to the integral equation (2.6), (2.7) respectively, it follows that
\begin{equation}\begin{split}\|E\|_{X_{S}^{k,b_1}}&\leq\|\psi_1(t)S(t)E_0\|_{X_{S}^{k,b_1}}+\|\psi_TS\ast_Rf_1\|_{X_{S}^{k,b_1}}
\\&\leq C\|E_0\|_{H^k}+CT^{1-b_1+b'}\|f_1\|_{X_{S}^{k,b'}},
\end{split}\end{equation}
and
\begin{equation}\begin{split}\|\varphi_{\pm}\|_{X_{W}^{l,b}}&\leq\|\psi_1(t)W(t)\varphi_{0\pm}\|_{X_{W}^{l,b}}+
\left\|\psi_T[W\ast_R(f+\Lambda^{-1}\psi_{2T}(\varphi_++\varphi_-))]\right\|_{X_{W}^{l,b}}
\\&\leq
C\|\varphi_{0\pm}\|_{H^l}+CT^{1-b-c}\|f\|_{X_{W}^{l,-c}}+CT\|\Lambda^{-1}\psi_{2T}(\varphi_++\varphi_-)\|_{X_{W}^{l,b}},
\end{split}\end{equation}
where $-\frac{1}{2}<b'\leq0\leq b_1\leq b'+1$ and $0\leq b\leq1-c$,
$c\in [0,\frac{1}{2}[$.
\par
 Let $b'=-c_1$, if $k,l,b$ satisfies the condition of Proposition 3.1 in \cite{GTV},
 then the second term of (2.14R) can be dealt with as follows
\begin{equation}\begin{split}T^{1-b_1+b'}\|f_1\|_{X_{S}^{k,b'}}&=T^{1-b_1-c_1}\|f_1\|_{X_{S}^{k,-c_1}}\\&\leq
CT^{1-b_1-c_1}\Big{(}T^{\gamma(b+b_1+c_1)}\|E\|_{X_{S}^{k,b_1}}\|\varphi_\pm\|_{X_{W}^{l,b}}\\&+\|\psi_{3T}^3(E\otimes
B)\|_{X_{S}^{k,-c_1}}\Big{)}.
\end{split}\end{equation}
Similarly, the second term of (2.15R) can be estimated by
\begin{equation}\begin{split}\|f\|_{X_{W}^{l,-c}}\leq
CT^{\gamma(c+2b_1)}\|E\|_{X_{S}^{k,b_1}}^2.
\end{split}\end{equation}
The solution $(E,\varphi_\pm)$ is locally well-posedness in the
space $X_{S}^{k,b_1}\times X_{W}^{l,b}$, if $\psi_{3T}(E\otimes B)$
satisfies \begin{equation} \|\psi_{3T}^3(E\otimes
B)\|_{X_{S}^{k,-c_1}}\leq C\|E\|_{X_{S}^{k,b_1}}^3.
\end{equation}
 If $b,b_1>\frac{1}{2}$, using Sobolev embedding theorem for time, then the
 solution satisfies
\begin{equation}(E,\varphi_\pm,\partial_t\varphi_\pm)\in\mathcal {C}([0,T[;H^k)\times\mathcal
{C}([0,T[;H^l)\times\mathcal {C}([0,T[;H^{l-1}).\end{equation} While
$b,\;b_1\leq\frac{1}{2}$, in order to obtain the continuity of time
(2.19), we need to prove
\begin{equation}\psi_{3T}^3(E\otimes B)\in
Y_S^k.\end{equation} In the next section, we devote to getting the
estimates (2.18), (2.10) and producing additional power of $T$ in
the process.
\section{The Nonlinear Estimates}
If the magnetic field $B$ satisfies the case $(A_1)$, then it
follows that
\begin{equation}B=\frac{i\eta}{(\Delta+\beta
I)}\nabla\times\nabla\times(E\otimes \overline{E}). \end{equation}
Substituting (3.1) into (2.18L), in order to derive (2.18R), it is
sufficient to show
\begin{equation}
|\mathcal {Q}|\leq
CT^{\vartheta_1}\|v\|_{L^2}\|v_1\|_{L^2}\|v_2\|_{L^2}\|v_3\|_{L^2},
\end{equation}
where$$\mathcal
{Q}=\int\frac{\langle\xi_1\rangle^k\hat{v}\hat{v}_1\hat{v}_2\hat{v}_3}{\langle\sigma\rangle^{b_1}\langle\sigma_1\rangle^{c_1}
\langle\widetilde{\sigma_2}\rangle^{b_1}\langle\sigma_3\rangle^{b_1}\langle\xi\rangle^k\langle
z_2\rangle^k\langle\xi_3\rangle^k},$$ $\xi=\xi_1-\xi_2$,
$z_2=\xi_2-\xi_3$, $\tau=\tau_1-\tau_2,s_2=\tau_2-\tau_3$,
$\widetilde{\sigma_2}=s_2+|z_2|^2$, $\sigma_i=\tau_i+|\xi|^2$,
$\hat{v}=\langle\xi\rangle^k\langle\sigma\rangle^{b_1}\hat{E}(\xi,\tau),\hat{v}_2=\langle
z_2\rangle^k\langle\widetilde{\sigma_2}\rangle^{b_1}\hat{E}(z_2,s_2),\hat{v}_3=\langle\xi_3\rangle^k\langle\sigma_3
\rangle^{b_1}\hat{E}(\xi_3,\tau_3)$ and $v_1\in L^2$. \par
Similarly, in order to estimate (2.20), we only to establish
\begin{equation}
|\mathcal {R}|\leq
CT^{\vartheta_2}\|v\|_{L^2}\|v_1\|_{L^2}\|v_2\|_{L^2}\|v_3\|_{L^2},
\end{equation}
where$$\mathcal
{R}=\int\frac{\langle\xi_1\rangle^k\hat{v}\hat{v}_1\hat{v}_2\hat{v}_3}{\langle\sigma\rangle^{b_1}\langle\sigma_1\rangle
\langle\widetilde{\sigma_2}\rangle^{b_1}\langle\sigma_3\rangle^{b_1}\langle\xi\rangle^k\langle
z_2\rangle^k\langle\xi_3\rangle^k},$$ with $v_1\in L^2_x$, the other
notation is the same as the above (3.2).
\par
On the other hand, if
$$A=-\gamma\frac{\partial}{\partial
t}\int_{\mathbb{R}^d}\frac{B(t,y)}{|x-y|^2}dy,$$ plugging it into
system $(1.1)_3$, after taking the partial Fourier transformation
with respect to the space variable, for $d>2$,
\begin{equation}|\xi|^2\widehat{B}(\xi)-i\eta\xi\times\xi\widehat{(E\otimes\overline{E})}(\xi)+\gamma
c_d|\xi|^{2-d}\widehat{B_t}(\xi)=0,
\end{equation}
where we have used the equality
$\widehat{|\cdot|^{-2}}=c_d|\cdot|^{2-d}$, for $d>2$. Consequently,
 \begin{equation}\partial_t\widehat{B}+\frac{1}{\gamma c_d}|\xi|^d\widehat{B}-\frac{i\eta}{\gamma
 c_d}|\xi|^{d-2}\xi\times\xi\widehat{(E\otimes\overline{E})}(\xi)=0.
\end{equation}
Therefore, we can solve $\widehat{B}$ as follows
$$\widehat{B}=e^{-\frac{1}{\gamma c_d}|\xi|^dt}\widehat{B}_0+\int_0^te^{-\frac{1}{\gamma c_d}|\xi|^d(t-\tau)}
\frac{i\eta}{\gamma
c_d}|\xi|^{d-2}\xi\times\xi\widehat{(E\otimes\overline{E})}d\tau,$$
which is equivalent to
\begin{equation}B=e^{-\frac{t}{\gamma c_d}(-\Delta)^{\frac{d}{2}}}B_0+\int_0^te^{-\frac{1}{\gamma c_d}(-\Delta)^{\frac{d}{2}}(t-\tau)}
\frac{i\eta}{\gamma
c_d}(-\Delta)^{\frac{d-2}{2}}(\nabla\times\nabla\times(E\otimes\overline{E}))d\tau.
\end{equation}
By virtue of the fractional parabolic equation theory, we only  to
estimate (3.2), (3.3). In order to prove (3.2), (3.3), we first give
the following lemma.
\begin{lemma3} Let $b_0>\frac{1}{2}$ and $\gamma\in[0,1]$. If
$a,a_1,a_2,a_3\geq0$  and $0<\eta,\eta_i\leq1,i=1,2,3$ satisfy
\begin{equation}(1-\gamma)[\eta
a+\eta_1a_1+\eta_2a_2+\eta_3a_3]=2b_0.
\end{equation}
Let $v,v_1,v_2,v_3\in L^2$ such that $\mathcal
{F}^{-1}(\langle\tau+|\xi|\rangle^{\gamma
a}\hat{v}),\mathcal{F}^{-1}(\langle\tau+|\xi_i|\rangle^{\gamma
a_i}\hat{v}_i),i=1,2$ have support in $|t|\leq CT$. Then the
following inequalities hold
\begin{equation}
\int\frac{\hat{v}\hat{v}_1\hat{v}_2\hat{v}_3}{\langle\sigma\rangle^{a}\langle\sigma_1\rangle^{a_1}
\langle\widetilde{\sigma_2}\rangle^{a_2}\langle\sigma_3\rangle^{a_3}\langle\xi\rangle^m}\leq
CT^{\theta}\|v\|_{L^2}\|v_1\|_{L^2}\|v_2\|_{L^2}\|v_3\|_{L^2},
\end{equation}
if $m\geq
d-(1-\gamma)[(1-\eta_1)a_1+(1-\eta_2)a_2+(1-\eta_3)a_3]/b_0\geq0$.
\begin{equation}
\int\frac{\hat{v}\hat{v}_1\hat{v}_2\hat{v}_3}{\langle\sigma\rangle^{a}\langle\sigma_1\rangle^{a_1}
\langle\widetilde{\sigma_2}\rangle^{a_2}\langle\sigma_3\rangle^{a_3}\langle\xi_3\rangle^{m_3}}\leq
CT^{\theta}\|v\|_{L^2}\|v_1\|_{L^2}\|v_2\|_{L^2}\|v_3\|_{L^2},
\end{equation}
if $m_3\geq
d-(1-\gamma)[(1-\eta)a+(1-\eta_1)a_1+(1-\eta_2)a_2]/b_0\geq0$.
\begin{equation}
\int\frac{\hat{v}\hat{v}_1\hat{v}_2\hat{v}_3}{\langle\sigma\rangle^{a}\langle\sigma_1\rangle^{a_1}
\langle\widetilde{\sigma_2}\rangle^{a_2}\langle\sigma_3\rangle^{a_3}\langle
\xi\rangle^m\langle\xi_3\rangle^{m_3}}\leq
CT^{\theta}\|v\|_{L^2}\|v_1\|_{L^2}\|v_2\|_{L^2}\|v_3\|_{L^2},
\end{equation}
 if $m+m_3\geq d-(1-\gamma)[(1-\eta_1)a_1+(1-\eta_2)a_2]/b_0\geq0$.
\begin{equation}
\int\frac{\hat{v}\hat{v}_1\hat{v}_2\hat{v}_3}{\langle\sigma\rangle^{a}\langle\sigma_1\rangle^{a_1}
\langle\widetilde{\sigma_2}\rangle^{a_2}\langle\sigma_3\rangle^{a_3}\langle
z_2\rangle^{m_2}\langle\xi_3\rangle^{m_3}}\leq
CT^{\theta}\|v\|_{L^2}\|v_1\|_{L^2}\|v_2\|_{L^2}\|v_3\|_{L^2},
\end{equation}
if $m_2+m_3\geq d-(1-\gamma)[(1-\eta)a+(1-\eta_1)a_1]/b_0\geq0$.
\begin{equation}
\int\frac{\hat{v}\hat{v}_1\hat{v}_2\hat{v}_3}{\langle\sigma\rangle^{a}\langle\sigma_1\rangle^{a_1}
\langle\widetilde{\sigma_2}\rangle^{a_2}\langle\sigma_3\rangle^{a_3}\langle
z_2\rangle^{m_2}\langle\xi\rangle^{m}}\leq
CT^{\theta}\|v\|_{L^2}\|v_1\|_{L^2}\|v_2\|_{L^2}\|v_3\|_{L^2},
\end{equation}
if $m+m_2\geq d-(1-\gamma)[(1-\eta_1)a_1+(1-\eta_3)a_3]/b_0\geq0$,
where $\theta=\gamma(a+a_1+a_2+a_3)$.
\end{lemma3}
\textit{Proof.} We first prove (3.8), taking advantage of H\"{o}lder
inequality in space time to (3.8L), we have
\begin{equation}\begin{split}
(3.8L)&\leq\|\mathcal
{F}^{-1}(\langle\xi\rangle^{-m}\langle\sigma\rangle^{-a}|\hat{v}|)\|_{L_t^q(L^r_x)}\|\mathcal
{F}^{-1}(\langle\widetilde{\sigma_2}\rangle^{-a_2}|\hat{v}_2|)\|_{L_t^{q_2}(L^{r_2}_x)}\\&\times
\prod_{i=1,3}\|\mathcal
{F}^{-1}(\langle\sigma_i\rangle^{-a_i}|\hat{v}_i|)\|_{L_t^{q_i}(L^{r_i}_x)}\\&
\leq C\|\mathcal
{F}^{-1}(\langle\xi\rangle^{-m}\langle\sigma\rangle^{-a}|\hat{v}|)\|_{L_t^q(L^r_x)}
\prod_{i=1}^3T^{\gamma a_i}\|v_i\|_{L^2},
\end{split}\end{equation}
where the coefficient satisfies
\begin{equation}\frac{1}{q}+\sum_{i=1}^3\frac{1}{q_i}=1,\quad\frac{1}{r}+\sum_{i=1}^3\frac{1}{r_i}=1,
\end{equation}
the second inequality comes from Lemma 3.1 in \cite{GTV} with the
coefficient satisfies for $i=1,2,3$
\begin{equation}\frac{2}{q_i}=1-\eta_i(1-\gamma)a_i/b_0,\quad\delta_i:=(\frac{d}{2}-\frac{d}{r_i})=(1-\eta_i)
(1-\gamma)a_i/b_0,
\end{equation}
Note that if $m\geq \delta(r):=\frac{d}{2}-\frac{d}{r}$, then
$H^m\hookrightarrow L^r$, consequently
\begin{equation}\begin{split}
\|\mathcal
{F}^{-1}(\langle\xi\rangle^{-m}\langle\sigma\rangle^{-a}\hat{v})\|_{L_t^q(L_x^r)}&\leq\|\mathcal
{F}^{-1}(\langle\xi\rangle^{-m}\langle\sigma\rangle^{-a}\hat{v})\|_{L_t^q(H_x^m)}\\&\leq\|\mathcal
{F}^{-1}(\langle\sigma\rangle^{-a}\hat{v})\|_{L_t^q(L_x^2)}\\&\leq
CT^{\gamma a}\|v\|_{L^2},
\end{split}\end{equation}
the last inequality comes from Lemma 2.4 in \cite{GTV} with
$\frac{2}{q}=1-(1-\gamma)\eta a/b_0$. Plugging (3.16) into (3.13),
one can easily get (3.8). Moreover, we have
$$(1-\gamma)[\eta a+\eta_1a_1+\eta_2a_2+\eta_3a_3]=2b_0$$
and$$m\geq\delta(r)=d-(\delta_1+\delta_2+\delta_3)=
d-(1-\gamma)[(1-\eta_1)a_1+(1-\eta_2)a_2+(1-\eta_3)a_3]/b_0\geq0.$$
Next, we will estimate (3.9) as follows
\begin{equation}\begin{split}
(3.9L)&\leq\|\mathcal
{F}^{-1}(\langle\sigma\rangle^{-a}|\hat{v}|)\|_{L_t^q(L^r_x)}
\|\mathcal
{F}^{-1}(\langle\sigma_1\rangle^{-a_1}|\hat{v}_1|)\|_{L_t^{q_1}(L^{r_1}_x)}\\&\times\|\mathcal
{F}^{-1}(\langle\widetilde{\sigma_2}\rangle^{-a_2}|\hat{v}_2|)\|_{L_t^{q_2}(L^{r_2}_x)}\|\mathcal
{F}^{-1}(\langle\sigma_3\rangle^{-a_3}\langle\xi_3\rangle^{-m_3}|\hat{v}_3|)\|_{L_t^{q_3}(L^{r_{3}}_x)}\\&
\leq C\|\mathcal
{F}^{-1}(\langle\sigma_3\rangle^{-a_3}\langle\xi_3\rangle^{-m_3}|\hat{v}_3|)\|_{L_t^{q_3}(L^{r_{3}}_x)}
T^{\gamma a}\|v\|_{L^2}\prod_{i=1}^2T^{\gamma
a_i}\|v_i\|_{L^2}\\&\leq
CT^{\gamma(a+a_1+a_2+a_3)}\|v\|_{L^2}\|v_1\|_{L^2}\|v_2\|_{L^2}\|v_3\|_{L^2},
\end{split}\end{equation}
where the constants satisfy (3.14), (3.15) and
$m_3\geq\delta_3\geq0$, which derives (3.7) and $$m_3\geq
d-(1-\gamma)[(1-\eta)a+(1-\eta_1)a_1+(1-\eta_2)a_2]/b_0\geq0.$$ We
now show (3.10), applying H\"{o}lder inequality in space time to
(3.10L) to yield
\begin{equation}\begin{split}
(3.10L)&\leq\|\mathcal
{F}^{-1}(\langle\xi\rangle^{-m}\langle\sigma\rangle^{-a}|\hat{v}|)\|_{L_t^q(L^r_x)}\|\mathcal
{F}^{-1}(\langle\xi_3\rangle^{-m_3}\langle\sigma_3\rangle^{-a_3}\\&|\hat{v}_3|)\|_{L_t^q(L^r_x)}\times
\prod_{i=1}^2\|\mathcal
{F}^{-1}(\langle\sigma_i\rangle^{-a_i}|\hat{v}_i|)\|_{L_t^{q_i}(L^{r_i}_x)}\\&\leq
CT^{\gamma(a+a_1+a_2+a_3)}\|v\|_{L^2}\|v_1\|_{L^2}\|v_2\|_{L^2}\|v_3\|_{L^2},
\end{split}\end{equation}
where we used (3.16) and \begin{equation}\begin{split} \|\mathcal
{F}^{-1}(\langle\xi_3\rangle^{-m_3}\langle\sigma_3\rangle^{-a_3}\hat{v}_3)&\|_{L_t^{q_3}(L_x^{r_3})}\\&\leq\|\mathcal
{F}^{-1}(\langle\xi_3\rangle^{-m_3}\langle\sigma_3\rangle^{-a_3}\hat{v}_3)\|_{L_t^{q_3}(H_x^{m_3})}\\&\leq\|\mathcal
{F}^{-1}(\langle\sigma_3\rangle^{-a_3})\|_{L_t^{q_3}(L_x^2)}\\&\leq
CT^{\gamma a_3}\|v\|_{L^2},
\end{split}\end{equation}
with $m_3\geq\delta_3:=\frac{d}{2}-\frac{d}{r_3}\geq0$,
$\frac{2}{q_3}=1-\frac{\eta_3(1-\gamma)a_3}{b_0}$. Thus, we have
$$m+m_3\geq\delta+\delta_3=d-(\delta_1+\delta_2)=d-\sum_{i=1}^2(1-\gamma)(1-\eta_i)a_i/b_0.$$
Similarly, one can easily check that (3.11) and (3.12). This
completes the proof of Lemma 3.1.\hspace{\fill}$\blacksquare$

\begin{lemma3} Assume $b_0>\frac{1}{2},\gamma\in[0,1]$ and $\eta,\eta_i\in(0,1]$,
$0<b_1,c_1<b_0$. Suppose the function $\mathcal
{F}^{-1}(\langle\sigma\rangle^{-b_1}\hat{v})$ and
$\mathcal{F}^{-1}(\langle\sigma_i\rangle^{-a_i}\hat{v}_i),i=1,2,3,$
have support in $t\leq CT$. If
$$(1-\gamma)[(\eta+\eta_2+\eta_3)b_1+\eta_1c_1]=2b_0, \qquad \text{and}$$
\begin{equation*}
\left\{\begin{array}{ll}2k\geq
d-(1-\gamma)[(1-\eta_1)c_1+(1-\eta_2)b_1]/b_0\geq0,
\\2k\geq d-(1-\gamma)[(1-\eta)b_1+(1-\eta_1)c_1]/b_0\geq0,\\
2k\geq d-(1-\gamma)[(1-\eta_1)c_1+(1-\eta_3)b_1]/b_0\geq0.
\end{array}\right.
\end{equation*}
Then the estimate (3.2) holds for all $T\leq T_0<\infty$ with
\begin{equation}|\mathcal {Q}|\leq
CT^{\gamma(3b_1+c_1)}\|v\|_{L^2}\|v_1\|_{L^2}\|v_2\|_{L^2}\|v_3\|_{L^2}.
\end{equation}
\end{lemma3}
\textit{Proof.}   In order to derive (3.20), we divide the
integration region into two subregions:
\par
Case 1, if $|\xi_1|\leq2|z_2|$, we estimate the contribution
$\mathcal {Q}_1$ of that region to $\mathcal {Q}$ by
\begin{equation}
|\mathcal{Q}_1|\leq\int\frac{|\hat{v}\hat{v}_1\hat{v}_2\hat{v}_3|}
{\langle\sigma\rangle^{b_1}\langle\sigma_1\rangle^{c_1}\langle\widetilde{\sigma_2}
\rangle^{b_1}\langle\sigma_3\rangle^{b_1}\langle\xi\rangle^k\langle\xi_3\rangle^k}.
\end{equation}
Thanks to (3.10) of Lemma 3.1 with
$(a,a_1,a_2,a_3,m,m_3)=(b_1,c_1,b_1,b_1,k,k)$, it follows that
\begin{equation}
|\mathcal {Q}_1|\leq
CT^{\gamma(3b_1+c_1)}\|v\|_{L^2}\|v_1\|_{L^2}\|v_2\|_{L^2}\|v_3\|_{L^2},
\end{equation}
and$$2k\geq d-(1-\gamma)[(1-\eta_1)c_1+(1-\eta_2)b_1]/b_0\geq0.$$
 Case 2,
if $|\xi_1|\geq2|z_2|$, then
$\frac{|\xi_1|}{2}\leq|\xi_1-z_2|\leq\frac{3}{2}|\xi_1|$, we
estimate the contribution $\mathcal {Q}_2$ of that region to
$\mathcal {Q}$ by
\begin{equation}\begin{split}
|\mathcal {Q}_2|&\leq
C\int\frac{\langle\xi_1-z_2\rangle^{k}|\hat{v}\hat{v}_1\hat{v}_2\hat{v}_3|}{\langle\sigma\rangle^{b_1}\langle\sigma_1\rangle^{c_1}
\langle\widetilde{\sigma_2}\rangle^{b_1}\langle\sigma_3\rangle^{b_1}\langle
\xi\rangle^k\langle z_2\rangle^k\langle\xi_3\rangle^k}\\&\leq
C\int\frac{|\hat{v}\hat{v}_1\hat{v}_2\hat{v}_3|}{\langle\sigma\rangle^{b_1}\langle\sigma_1\rangle^{c_1}
\langle\widetilde{\sigma_2}\rangle^{b_1}\langle\sigma_3\rangle^{b_1}\langle
z_2\rangle^k\langle\xi_3\rangle^k}\\&+C\int\frac{|\hat{v}\hat{v}_1\hat{v}_2\hat{v}_3|}{\langle\sigma\rangle^{b_1}
\langle\sigma_1\rangle^{c_1}\langle\widetilde{\sigma_2}\rangle^{b_1}\langle\sigma_3\rangle^{b_1}\langle
z_2\rangle^k\langle\xi\rangle^k}\\&:=\mathcal {Q}_{21}+\mathcal
{Q}_{22}.
\end{split}\end{equation}
By virtue of  (3.11) of Lemma 3.1 with
$(a,a_1,a_2,a_3,m,m_3)=(b_1,c_1,b_1,b_1,k,k)$, we deduce that
\begin{equation}
|\mathcal {Q}_{21}|\leq
CT^{\gamma(3b_1+c_1)}\|v\|_{L^2}\|v_1\|_{L^2}\|v_2\|_{L^2}\|v_3\|_{L^2},
\end{equation}
and$$2k\geq d-(1-\gamma)[(1-\eta)b_1+(1-\eta_1)c_1]/b_0\geq0.$$
Similarly, we can estimate $\mathcal {Q}_{22}$ as follows
\begin{equation}
|\mathcal {Q}_{22}|\leq
CT^{\gamma(3b_1+c_1)}\|v\|_{L^2}\|v_1\|_{L^2}\|v_2\|_{L^2}\|v_3\|_{L^2},
\end{equation}
and$$2k\geq d-(1-\gamma)[(1-\eta_1)c_1+(1-\eta_3)b_1]/b_0\geq0.$$
This completes the proof of Lemma 3.2.\hspace{\fill}$\blacksquare$
\begin{lemma3} Let $b_0>\frac{1}{2}$ and $\gamma\in[0,1]$. Given $\eta,\eta_i\in(0,1]$,
 and $0<b_1,c_1<b_0$. Assume that the function $\mathcal
{F}^{-1}(\langle\sigma\rangle^{-b_1}\hat{v})$ and
$\mathcal{F}^{-1}(\langle\sigma_i\rangle^{-a_i}\hat{v}_i),i=1,2,3,$
have support in $t\leq CT$. If
$$(1-\gamma)[(\eta+\eta_2+\eta_3)b_1+\eta_1]=2b_0, \qquad \text{and}$$
\begin{equation*}
\left\{\begin{array}{ll}2k\geq
d-(1-\gamma)[(1-\eta_1)+(1-\eta_2)b_1]/b_0\geq0,
\\k\geq d-(1-\gamma)[(2-\eta-\eta_2)b_1+1-\eta_1]/b_0\geq0,\\
k\geq d-(1-\gamma)[1-\eta_1+(2-\eta_2-\eta_3)b_1]/b_0\geq0.
\end{array}\right.
\end{equation*}
Then the estimate (3.3) holds for all $T\leq T_0<\infty$ with
\begin{equation}|\mathcal {R}|\leq
CT^{\gamma(3b_1+1)}\|v\|_{L^2}\|v_1\|_{L^2}\|v_2\|_{L^2}\|v_3\|_{L^2}.
\end{equation}
\end{lemma3}
\textit{Proof.} In order to estimate (3.26), we divide the
integration region into two subregions:
\par Region
$|\xi_1|\leq2|z_2|$: We estimate the contribution $\mathcal {R}_1$
of that region to $\mathcal {R}$ by
\begin{equation}
|\mathcal{R}_1|\leq\int\frac{|\hat{v}\hat{v}_1\hat{v}_2\hat{v}_3|}
{\langle\sigma\rangle^{b_1}\langle\sigma_1\rangle\langle\widetilde{\sigma_2}
\rangle^{b_1}\langle\sigma_3\rangle^{b_1}\langle\xi\rangle^k\langle\xi_3\rangle^k}.
\end{equation}
Thanks to (3.10) of Lemma 3.1 with
$(a,a_1,a_2,a_3,m,m_3)=(b_1,1,b_1,b_1,k,k)$, one can easily check
that
\begin{equation}
|\mathcal {R}_1|\leq
CT^{\gamma(3b_1+1)}\|v\|_{L^2}\|v_1\|_{L^2}\|v_2\|_{L^2}\|v_3\|_{L^2},
\end{equation}
and$$2k\geq d-(1-\gamma)[(1-\eta_1)+(1-\eta_2)b_1]/b_0\geq0.$$
Region $|\xi_1|\geq2|z_2|$: Since
$\frac{|\xi_1|}{2}\leq|\xi_1-z_2|\leq\frac{3}{2}|\xi_1|$, we
estimate the contribution $\mathcal {R}_2$ of that region to
$\mathcal {R}$ by
\begin{equation}\begin{split}
|\mathcal {R}_2|&\leq
C\int\frac{|\hat{v}\hat{v}_1\hat{v}_2\hat{v}_3|}{\langle\sigma\rangle^{b_1}\langle\sigma_1\rangle
\langle\widetilde{\sigma_2}\rangle^{b_1}\langle\sigma_3\rangle^{b_1}\langle
\langle\xi_3\rangle^k}\\&+C\int\frac{|\hat{v}\hat{v}_1\hat{v}_2\hat{v}_3|}{\langle\sigma\rangle^{b_1}
\langle\sigma_1\rangle\langle\widetilde{\sigma_2}\rangle^{b_1}\langle\sigma_3\rangle^{b_1}\langle\xi\rangle^k}\\&:=\mathcal
{R}_{21}+\mathcal {R}_{22}.
\end{split}\end{equation}
By virtue of (3.8), (3.9) of Lemma 3.1 with
$(a,a_1,a_2,a_3,m_3)=(b_1,1,b_1,b_1,k)$,
$(a,a_1,a_2,a_3,m)=(b_1,1,b_1,b_1,k)$, we can end up with
\begin{equation}
|\mathcal {R}_{21}|+|\mathcal {R}_{22}|\leq
CT^{\gamma(3b_1+1)}\|v\|_{L^2}\|v_1\|_{L^2}\|v_2\|_{L^2}\|v_3\|_{L^2},
\end{equation}
with the efficient satisfies $k\geq
d-(1-\gamma)[(2-\eta-\eta_2)b_1+1-\eta_1]/b_0\geq0,$ and $$k\geq
d-(1-\gamma)[1-\eta_1+(2-\eta_2-\eta_3)b_1]/b_0\geq0.$$ This
completes the proof of Lemma 3.3.\hspace{\fill}$\blacksquare$
\par
Analogous to Proposition 1.1 in \cite{GTV}, in view of Lemma 3.2,
Lemma 3.3, we have the  following well-posedness result.
\begin{theorem3} Let the space dimension $d>1$. Assume $k,l$ satisfy
\begin{equation}
k\in[l,l+1],\;l>\frac{d}{2}-2,\;2k-(l+1)>\frac{d}{2}-2,\qquad
\text{if}\quad d\geq4,\end{equation}
\begin{equation}
k\in[l,l+1],\;\qquad l\geq0,\;\quad2k-(l+1)\geq0,\qquad
\text{if}\quad d=2,3.\end{equation} Then GZ system (1.1) in the case
 $(A_1)$ (or in the case $(A_2)$ if $d>2$) with the initial data $(E_0,n_0,\partial_tn_0)\in
 H^k\times
H^l\times H^{l-1}$ is locally  well-posedness in $X_S^{k,b_1}\times
X_W^{l,b}\times X_W^{l-1,b}$ with the $B_0\in H^k$ for suitable
$b_1,b$ close to $\frac{1}{2}$. Moreover, the solutions satisfy
\begin{equation}
(E,n,\partial_tn)\in \mathcal {C}([0,T[;H^k\times H^l\times
H^{l-1}).
\end{equation}
\end{theorem3}
\textit{Proof.} By virtue of Lemma 3.2 and Lemma 3.3. Substituting
(2.18) into (2.16), then plugging (2.16), (2.17) into (2.14), (2.15)
respectively, it follows from (2.14) and (2.15) that
\begin{equation}\begin{split}
\|E\|_{X_S^{k,b_1}}&\leq
C\|E_0\|_{H^k}+CT^{1-b_1-c_1}\Big{(}T^{\gamma_1(b+b_1+c_1)}\|\psi_{2T}E\|_{X_S^{k,b_1}}\\&\qquad\times\|\psi_{2T}\varphi_{\pm}\|_{X_W^{l,b}}
+T^{\gamma_2(3b_1+c_1)}\|\psi_{3T}E\|_{X_S^{k,b_1}}^3\Big{)},
\end{split}
\end{equation}
and
\begin{equation}\begin{split}
\|\varphi_\pm\|_{X_W^{l,b}}\leq
C\|\varphi_{0\pm}\|_{H^l}&+CT^{1-b-c}\Big{(}T^{\gamma_3(2b_1+c)}\\&\times\|\psi_{2T}E\|_{X_S^{k,b_1}}^2+T^{b+c}\|\psi_{2T}\varphi_{\pm}\|_{X_W^{l,b}}
\Big{)},
\end{split}
\end{equation}
where $k,l$ satisfy (3.31), (3.32) and the condition of Lemma 3.2.
Consequently, we solve (3.34) and (3.35) by the contraction mapping
argument for small enough time $T$ in space $X_S^{k,b_1}\times
X_W^{l,b}$. It remains only to be proved that under condition (3.31)
and (3.32), we can choose $b_1,c_1$, $\eta_i\in(0,1]$ and
$\gamma\in[0,1],i=1,2,3$ satisfying the assumptions of Lemma 3.2 and
Lemma 3.3 if needed. Since (3.31), (3.32) and the assumptions of
Lemma 3.2 and Lemma 3.3 are consistency condition, which it is not
difficult to check.
\par
At this point we have obtained the existence, uniqueness and
continuity of local solution in time for the cut-off equation (2.6)
and (2.7). Similar to the method on page 415--416 in \cite{GTV}, one
can easily check that the solutions are in fact not depend with the
cut-off time. This completes the proof of Theorem
3.1.\hspace{\fill}$\blacksquare$
\begin{remark3}
Although the form of system (1.1) is more complex than it of
Eq.(1.2), the result well-posedness of Theorem 3.1 to the GZ system
(1.1) with magnetic field $(A_1)$ (or $(A_2)$, if $d>2$) is the same
as the Proposition 1.1 \cite{GTV}, which is proved by J. Ginibre, Y.
Tsutsumi and G. Velo in 1996.
\end{remark3}
\begin{remark3}
If the space dimension $d=1$, then $E\otimes B=0$, the system (1.1)
with the case $(A_1)$ becomes Eq.(1.2). Thus the GZ system (1.1) is
locally well-posedness for $(E_0,n_0,\partial_tn_0)\in H^k\times
H^l\times H^{l-1}$, if the indexes $k$ and $l$ satisfy
$$-\frac{1}{2}\leq k-l\leq1,\quad 0\leq l+\frac{1}{2}\leq2k.$$
Moreover, if the initial data $(E_0,n_0,\partial_tn_0)\in H^1\times
L^2\times H^{-1}$, then there exists a global solution
$(E,n,\partial_tn)$, which satisfies
$$(E,n,\partial_tn)\in \mathcal {C}(\mathbb{R}^+;H^1)\times\mathcal {C}
(\mathbb{R}^+;L^2)\times\mathcal {C}(\mathbb{R}^+;H^{-1}).$$
\end{remark3}
\begin{remark3} The initial data $B_0\in H^k$ in Theorem 3.1
is necessary. In fact, if $A$ satisfies the case $(A_1)$, then we
can derive from (3.6)
that$$\|B\|_{X_S^{k,b_1}}\leq\|B_0\|_{H^k}+C\|E\|_{X_S^{k,b_1}}^2.$$
\end{remark3}
\begin{remark3}
By the energy estimation, C. Laurey \cite{La} proves the local
existence of solution in the spaces $$(E,n,B)\in \mathcal
{C}([0,T[;H^{s+1})\times\mathcal {C}([0,T[;H^{s})\times\mathcal
{C}([0,T[;H^{s+1}),$$ to system (1.1) with the case $(A_1)$,
$(A_2)$, for $d=2,3$. In fact, similarly, we also can prove for some
$T>0$ that the solutions satisfy
$$(E,n,B)\in \mathcal {C}([0,T[;B_{p,r}^{s+1})\times\mathcal {C}([0,T[;B_{p,r}^{s})\times\mathcal {C}([0,T[;B_{p,r}^{s+1}),$$
for $s>\frac{d}{2}$.
\end{remark3}

\section{The lower bound for the blow-up rate of blow-up solutions}
\par
In 1994, L. Glangetas and F. Merle \cite{GM} proved the following
form of self-similar blow-up solutions to equation (1.2) in
$\mathbb{R}^2$, i.e. \begin{equation}
\left\{\begin{array}{ll}E(t,x)=\frac{\omega}{T-t}e^{i(\theta+\frac{|x|^2-4\omega^2}{4(-T+t)})}
P\left(\frac{x\omega}{T-t}\right),
\\n(t,x)=\left(\frac{\omega}{T-t}\right)^2N\left(\frac{x\omega}{(T-t)}\right).\end{array}\right.
\end{equation}
where $\omega>\omega_0,\theta\in\mathbb{R},$
$P(x)=P(|x|),N(x)=N(|x|)$, and $(P,N)$ satisfies the elliptic
equation \begin{equation*} \left\{\begin{array}{ll}\Delta P-P=NP,\\
\frac{1}{(c_0\omega)^2}(r^2N_{rr}+6rN_r+6N)-\Delta N=\Delta
P^2\end{array}\right.
\end{equation*}
with $r=|x|, \Delta w=w_{rr}+\frac{1}{r}w_r$. \par The situation for
Zakharov equation in $\mathbb{R}^3$ is more complex. Until now,
there are no known explicit blow-up solutions, M. Landman, etc.
observed an asymptotic self-similar blow-up solution for Zakharov
equation (1.2) in $\mathbb{R}^3$ of the form \cite{LPSSW}
\begin{equation}
\left\{\begin{array}{ll}E(t,x)=\frac{2}{3(T-t)}e^{i(T-t)^{-1/3}}
P\left(\frac{|x|}{\sqrt{3}(T-t)^{2/3}}\right),
\\n(t,x)=\frac{1}{3(T-t)^{4/3}}N\left(\frac{|x|}{\sqrt{3}(T-t)^{2/3}}\right).\end{array}\right.
\end{equation}
where $P(x)=P(|x|),N(x)=N(|x|)$, and $(P,N)$ satisfies the elliptic
equation \begin{equation*} \left\{\begin{array}{ll}\Delta P-P=NP,\\
\frac{1}{(2}(2r^2N_{rr}+13rN_r+14N)=\Delta P^2.\end{array}\right.
\end{equation*}
\par
In this subsection, we consider the singular solution of the system
(1.1) in the case $(A_1)$ and $(A_2)$ in finite time, we will
establish the lower bound for the blow-up rate of the blow-up
solution to system (1.1).
\begin{theorem4} Let $k,l$ satisfy (3.31) and (3.32). Assume that the initial data $(E_0,B_0,n_0,n_1)$
belongs to $H^k\times H^k\times H^l\times H^{l-1}$. Then there
exists a time $T>0$ depending only on $(E_0,B_0,n_0,n_1)$ and a
unique solution $(E,n,\partial_tn)$ to system (1.1) with the initial
data $(E_0,n_0,n_1)$, which is guaranteed by Theorem 3.1. If the
solution blows up in finite time
 $T^*$ in the space
 $H^{k}\times H^{l}\times H^{l-1}$, then we have the lower bound for the
 blow-up rate of blow-up solution satisfies for any $\epsilon>0$
$$\left(\|E(t)\|_{H^{k}}+\|n(t)\|_{H^l}+\|\partial_tn(t)\|_{H^{l-1}}\right)>C\frac{1}{(T^*-t)^{\frac{3}{4}\vartheta_l-\epsilon}},$$
where $\vartheta_l=\frac{1}{4}(2l+4-d)$.
\end{theorem4}
 \textit{Proof.} Let $b=b_1,c=c_1$ and $\gamma_1=\gamma_3$, in view of (3.34) and (3.35), we have
\begin{equation}\begin{split}
\|E\|_{X_S^{k,b}}&\leq
C\|E_0\|_{H^k}+CT^{1-b-c+\gamma_1(2b+c)}\|\psi_{2T}E\|_{X_S^{k,b}}\\&\qquad\times\|\psi_{2T}\varphi_{\pm}\|_{X_W^{l,b}}
+CT^{1-b-c+\gamma_2(3b+c)}\|\psi_{3T}E\|_{X_S^{k,b}}^3\\& \leq
C\|E_0\|_{H^k}+CT^{2-3b-c+\gamma_1(2b+c)}\|E\|_{X_S^{k,b}}\|\varphi_{\pm}\|_{X_W^{l,b}}\\&\qquad\qquad+
CT^{\frac{5}{2}-4b-c+\gamma_2(3b+c)}\|E\|_{X_S^{k,b}}^3,
\end{split}
\end{equation}
and
\begin{equation}\begin{split}
\|\varphi_\pm\|_{X_W^{l,b}}&\leq
C\|\varphi_{0\pm}\|_{H^l}+CT^{1-b-c}\Big{(}T^{\gamma_1(2b+c)}\|\psi_{2T}E\|_{X_S^{k,b}}^2+T^{b+c}\|\psi_{2T}\varphi_{\pm}
\|_{X_W^{l,b}}\Big{)}\\&\leq
C\|\varphi_{0\pm}\|_{H^l}+C\Big{(}T^{2-3b-c+\gamma_1(2b+c)}\|E\|_{X_S^{k,b}}^2+T^{\frac{3}{2}-b}\|\varphi_{\pm}
\|_{X_W^{l,b}}\Big{)},
\end{split}
\end{equation}
the last inequality of (4.3), (4.4) comes from Remark 2.1 with
$q=2$, where $k,l$ satisfy (3.31), (3.32) and the other
constants satisfy \begin{equation}b_0>\frac{1}{2},\;0<b,c\leq
b_0,\;\gamma_i\in[0,1],\end{equation}
\begin{equation}(1-\gamma_1)\text{max}\{b,c\}\leq
b_0\leq(1-\gamma_1)(2b+c),\end{equation}
\begin{equation}l\geq\frac{d}{2}+1-(1-\gamma_1)(2b+c)/b_0,\end{equation} and
\begin{equation}
\left\{\begin{array}{ll}(1-\gamma_2)[(\eta+\eta_2+\eta_3)b_1+\eta_1c]=2b_0,\quad
\eta_i\in(0,1],
\\2k\geq
d-(1-\gamma_2)[(1-\eta_1)c+(1-\eta_2)b]/b_0\geq0,
\\2k\geq d-(1-\gamma_2)[(1-\eta)b+(1-\eta_1)c]/b_0\geq0,\\
2k\geq
d-(1-\gamma_2)[(1-\eta_1)c+(1-\eta_3)b]/b_0\geq0,\end{array}\right.
\end{equation}
In view of (4.7), in order to let $\gamma_1$ large enough, we choose
\begin{equation}
\gamma_1=1-\frac{(\frac{d}{2}+1-l)b_0}{2b+c}.
\end{equation}
Due to $c\in(0,\frac{1}{2}),b>\frac{1}{2}$. Substituting (4.9) into
(4.6) to yield
\begin{equation}\frac{d}{2}-2<l<\frac{d+2}{2}.\end{equation}
By virtue of (4.9) and (4.10), let
$b=\frac{1}{2}+\epsilon,b_0=\frac{1}{2}+\epsilon_0,0<\epsilon\leq\epsilon_0$,
it follows that
\begin{equation}\begin{split}2-3b-c+\gamma_1(2b+c)&=2-b+(l-1-\frac{d}{2})b_0\\&=\frac{4-d}{4}+\frac{l}{2}-
\epsilon-\left(\frac{d+2}{2}-l\right)\epsilon_0\\&:=\vartheta_l^-\\&<\frac{1}{4}(2l+4-d):=\vartheta_l.
\end{split}\end{equation}
Let $\frac{5}{2}-4b-c+\gamma_2(3b+c)=\frac{3}{2}\vartheta_l^-$,
which can be guaranteed by (4.8). Combining (4.3) with (4.4), in
view of H\"{o}lder inequality to give by
\begin{equation}\begin{split}
\|E\|_{X_S^{k,b}}&+\|\varphi_{\pm}\|_{X_W^{l,b}}\leq C(
\|E_0\|_{X_S^{k,b}}+\|\varphi_{0\pm}\|_{X_W^{l,b}})+CT^{\frac{3}{2}-b}\|\varphi_{\pm}\|_{X_W^{l,b}}\\&+CT^{\vartheta_l^-}(\|E\|_{X_S^{k,b}}+\|\varphi_{\pm}\|_
{X_W^{l,b}})^2+CT^{\frac{3}{2}\vartheta_l^-}(\|E\|_{X_S^{k,b}}+\|\varphi_{\pm}\|_
{X_W^{l,b}})^3\\&\leq
C(\|E_0\|_{X_S^{k,b}}+\|\varphi_{0\pm}\|_{X_W^{l,b}})+CT^{\frac{3}{2}\vartheta_l^-}(\|E\|_{X_S^{k,b}}+\|\varphi_{\pm}\|_
{X_W^{l,b}})^3\\&
\qquad\qquad\qquad\qquad+CT^{\frac{3}{2}-b}(\|E\|_{X_S^{k,b}}+\|\varphi_{\pm}\|_{X_W^{l,b}}).
\end{split}\end{equation}
Note that $\frac{3}{2}-b>0$, if $T$ is small enough such that
$CT^{\frac{3}{2}-b}\leq \frac{1}{2}$, then we have
\begin{equation}
\|E\|_{X_S^{k,b}}+\|\varphi_{\pm}\|_{X_W^{l,b}}\leq
C(\|E_0\|_{X_S^{k,b}}+\|\varphi_{0\pm}\|_{X_W^{l,b}})+CT^{\frac{3}{2}\vartheta_l^-}(\|E\|_{X_S^{k,b}}+\|\varphi_{\pm}\|_
{X_W^{l,b}})^3.
\end{equation}
 Next, we will infer a lower bound on the blow-up rate of blow-up solution. Denote by
$T^*$ the supremum of the existence time $T>0$ for which there
exists a solution $(E,n)$ of the Zakharov system (1.1) satisfying
$$(\|E\|_{X_S^{k,b}}+\|n\|_{X_W^{l,b}}+\|\partial_tn\|_{X_W^{l-1,b}})<\infty.$$
Then for all time $t\in [0,T)[$, the solutions satisfy
$$\|E(t)\|_{H^k}+\|n(t)\|_{
H^l}+\|\partial_tn(t)\|_{H^{l-1}}<\infty,$$ which is guaranteed by
the local well-posedness of Theorem 3.1.  By the maximality of
$T^*$, it follows that
$$\|E(t)\|_{L^\infty_{T^*}(H^k)}+\|n(t)\|_{L^\infty_{T^*}(
H^l)}+\|\partial_tn(t)\|_{L^\infty_{T^*}(H^{l-1})}=\infty.$$
Otherwise, the Cauchy problem of system (1.1) at time $T^*$ with the
initial data $(E(T^*,\cdot),n(T^*,\cdot))$ would be well-defined and
the local existence theory would extend the solution $(E,n)$ beyond
$T^*$. Thus, if $T^*<\infty$, the solution blows up and
$$\|E(t)\|_{H^k}+\|n(t)\|_{
H^l}+\|\partial_tn(t)\|_{H^{l-1}}\longrightarrow\infty\qquad
t\rightarrow T^*.$$
\par
Consider the solution $(E,n)$ posed at some time $t\in[0,T^*[$.
Assume for some $M$ such that
$$C(|E(t)\|_{H^k}+\|n(t)\|_{
H^l}+\|\partial_tn(t)\|_{H^{l-1}})+C(T-t)^{\frac{3}{2}\vartheta_l^-}M^3\leq
M.$$ Then $T<T^*$. Consequently, $\forall M>0$
$$C(\|E(t)\|_{H^k}+\|n(t)\|_{
H^l}+\|\partial_tn(t)\|_{H^{l-1}})+C(T^*-t)^{\frac{3}{2}\vartheta_l^-}M^3>
M.$$ Choosing $M=2C(\|E(t)\|_{H^k}+\|n(t)\|_{
H^l}+\|\partial_tn(t)\|_{H^{l-1}})$, we deduce that
$$C(T^*-t)^{\frac{3}{2}\vartheta_l^-}M^3>M,$$
which is equivalent to
$$(\\|E(t)\|_{H^k}+\|n(t)\|_{
H^l}+\|\partial_tn(t)\|_{H^{l-1}})>
C\frac{1}{(T^*-t)^{\frac{3}{4}\vartheta_l^-}}.$$ This completes the
proof of Theorem 4.1.\hspace{\fill}$\blacksquare$
\begin{corollary4}Under the assumption of Theorem 4.1. If we
neglect the magnetic $B$, then the classical Zakharov Eq.(1.2) is
locally well-posedness. If the solution $(E,n)$ blows up in finite
time $T^*$ in the space
 $H^{k}\times H^{l}\times H^{l-1}$, then we have the lower bound for the
 blow-up rate of blow-up solution satisfies for any $\epsilon>0$
$$\left(\|E\|_{H^{k}}+\|n\|_{H^l}+\|\partial_tn\|_{H^{l-1}}\right)>C\frac{1}{(T^*-t)^{\vartheta_l-\epsilon}},$$
where $\vartheta_l=\frac{1}{4}(2l+4-d)$.
\end{corollary4}
\textit{Proof.} As the process of (4.12), we have
\begin{equation}\begin{split}
\|E\|_{X_S^{k,b}}+\|\varphi_{\pm}\|_{X_W^{l,b}}&\leq
C(\|E_0\|_{H^k}+\|\varphi_{0\pm}\|_{H^l})+CT^{\vartheta_l^-}(\|E\|_{X_S^{k,b}}+\|\varphi_{\pm}\|_{X_W^{l,b}})^2\\&
+CT^{\frac{5}{2}-2b-c}(\|E\|_{X_S^{k,b}}+\|\varphi_{\pm}\|_{X_W^{l,b}}),
\end{split}\end{equation}
where $\vartheta_l=\frac{1}{4}(2l+4-d)$. Similarly, one can easily
get the lower bound for blow-up rate of blow-up solution to Eq.(1.2)
\begin{equation}
\left(\|E(t)\|_{H^{k}}+\|n(t)\|_{H^l}+\|\partial_tn(t)\|_{H^{l-1}}\right)>C\frac{1}{(T^*-t)^{\vartheta_l-\epsilon}}.
\end{equation}
This concludes the proof of Corollary
4.1.\hspace{\fill}$\blacksquare$

\begin{remark4} If we consider the self-similar blow-up solution
$(E, n)$ of (4.1) to Eq.(1.2) which blows up in a finite time $T^*$
in $\mathbb{R}^2$, then we obtain the blow-up rate of blow-up
solution $n$ satisfying
$$\|n\|_{\dot{H}^l}=\left(\frac{\omega}{T^*-t}\right)^{l+1}\|N\|_{\dot{H}^l}.$$
In \cite{Me1}, F. Merle prove the optimal lower bound of the blow-up
rate of the solution (E,n)in space $H^1\times L^2$ in 2D is
$C\frac{1}{(T^*-t)}$. However, for $d=3$, the homogeneous norm of
$n$ the asymptotic self-similar blow-up solution (4.2) is
$$\|n\|_{\dot{H}^l}=\frac{C}{(T^*-t)^{\frac{1}{3}(2l+1)}}\|N\|_{\dot{H}^l}.$$
\end{remark4}
\begin{remark4} As $d=3$, the result of Corollary 4.1 was obtained in \cite{CCS}.
If we consider the following asymptotic self-similar blow-up
solution to Eq.(1.2) with $c_0=1$ in 3D
\begin{equation}
\left\{\begin{array}{ll}E(t,x)=\frac{1}{T-t}
P\left(\frac{|x|}{(T-t)^{1/2}}\right)+i\frac{1}{T-t}Q\left(\frac{|x|}{(T-t)^{1/2}}\right),
\\n(t,x)=\frac{1}{T-t}N\left(\frac{|x|}{(T-t)^{1/2}}\right).\end{array}\right.
\end{equation}
where $P(x)=P(|x|),Q(x)=Q(|x|),N(x)=N(|x|)$, and $(P,N)$ satisfies
the ODEs
\begin{equation*} \left\{\begin{array}{ll}\Delta P+\frac{1}{2}rQ_r+Q=NP,\\
\Delta Q+\frac{1}{2}rP_r+P=NQ,\\
\frac{1}{4}r^2N_{rr}+\frac{7}{4}rN_r+2N)=\Delta
(P^2)+\Delta(Q^2).\end{array}\right.
\end{equation*}
The homogeneous norm of the solution (4.16) satisfies
$$\|n\|_{\dot{H}^l}=(\frac{1}{T-t})^{\frac{1}{4}(2l+1)}\|N\|_{\dot{H}^l},$$
and
$$\|E\|_{\dot{H}^k}=(\frac{1}{T-t})^{\frac{1}{4}(2k+1)}(\|P\|_{\dot{H}^k}+\|Q\|_{\dot{H}^k}).$$
The lower bound of the solution in Corollary 4.1 is almost up to the  optimal bound of
the asymptotic blow-up rate $\frac{1}{4}(2l+1)$. Until now,  we do
not find explicit blow-up solution to Eq.(1.2) in $\mathbb{R}^3$,
the blow-up rate of blow-up solution is open problem.
\end{remark4}
\section{The global existence of solution}
In this subsection, by the local well-posedness and conservation laws,  in the space dimension $d=2,3,4$, we shall establish the global solution of the GZ system with magnetic field in the case $(A_1)$, the results are
\begin{theorem5} Assume the initial data $(E_0,n_0,n_1)$ belong to the Sobolev space $H^k(\mathbb{R}^d)\times H^{k-1}(\mathbb{R}^d)\times H^{k-2}(\mathbb{R}^d),k\geq1,\;d=2,3$. Let
the initial data satisfy
\begin{equation}
\left\{\begin{array}{ll}(1+\frac{\eta}{2})K^4(2)\|E_0\|_{L^2}^2<1,&\text{if} \;\;d=2,\\
\|E_0\|_{L^2}^2 H_2(0)(1+\frac{\eta}{2})^2<\frac{4}{27K^8(3)},\\
 \text{and}\;\|\nabla E_0\|_{L^2}^2\leq H_2(0),&\text{if} \;\;d=3.
\end{array}\right.
\end{equation}
Then there exists a unique and global solution $$(E,n,\partial_tn)\in \mathcal{C}(\mathbb{R}^+;H^k(\mathbb{R}^2)\times H^{k-1}(\mathbb{R}^2)\times H^{k-2}(\mathbb{R}^2))$$ and
$$(E,n,\partial_tn)\in \mathcal{C}(\mathbb{R}^+;H^1(\mathbb{R}^3)\times L^2(\mathbb{R}^3)\times H^{-1}(\mathbb{R}^3))$$ to system (1.1) in the case $(A_1)$ with the initial data $(E_0,n_0,n_1)$. Moreover, if $k=1$, then the global solutions satisfy
\begin{equation}
\|E\|_{H^1}+\|n\|_{L^2}+\|V\|_{L^2}+\|\partial_tn\|_{H^{-1}}\leq C,
\end{equation}
uniformly bound for $t\in\mathbb{R}^+$.
If $k\geq2$ and $d=2$, the global solutions grow at most $k$-exponential bounds and satisfy
\begin{equation}
\|E\|_{H^k}+\|n\|_{H^{k-1}}+\|\partial_tE\|_{H^{k-2}}+\|\partial_tn\|_{H^{k-2}}\leq Ce^{e^{\cdots^{e^{ct}}}},
\end{equation}
where $c$ and $C$ are positive constants.
\end{theorem5}
\par
At first, in order to present the proof the Theorem 5.1, we recall the following two lemmas.
\begin{lemma5}\cite{We} Let the function $u\in H^1(\mathbb{R}^d),2\leq d<4$. Then we have
\begin{equation*}
\|u\|_{L^4}^4\leq K^4(d)\|u\|_{L^2}^{4-d}\|\nabla u\|_{L^2}^d,
\end{equation*}
where $K^4(d)=\frac{2}{\|\psi\|_{L^2}^2}$, the function $\psi$ is the ground state solution of
\begin{equation*}
\frac{d}{2}\Delta\psi+\frac{d-4}{2}\psi+\psi^3=0.
\end{equation*}
\end{lemma5}
\begin{lemma5} Given $f(t)$ be positive and continuous function on $\mathbb{R}^+$. Let $c_1,c_2>0$ and $k>1$ such that
$$f(t)\leq c_1+c_2f^k(t).$$ If the constants  $c_1,c_2>0$ and $k>1$ satisfy
$$c_1^{k-1}c_2<\frac{(k-1)^{k-1}}{k^k},\quad\text{and}\;f(0)\leq c_1,$$ then the function $f$ is uniformly bounded on $\mathbb{R}^+$.
\end{lemma5}
\par The proof of Lemma 5.2 is simple, which can be found in many books, we omit it here.\\
\textit{Proof of Theorem 5.1.} The local well-posedness of solution to system (1.1) is guaranteed by Theorem 3.1. In view of the conservation law (1.5), (1.6) and Lemma 5.1, without loss of generality, let $c_0=1$, one can easily check that
\begin{equation}\begin{split}
\|\nabla E\|_{L^2}^2+\frac{1}{4}\|n\|_{L^2}^2+\frac{1}{2}\|V\|_{L^2}^2\leq H_2(0)+(1+\frac{\eta}{2})K^{4}(d)\|E_0\|_{L^2}^{4-d}\|\nabla E\|_{L^2}^d,
\end{split}\end{equation}
where we have used Young's inequality, H\"{o}lder's inequality and
\begin{equation*}\begin{split}\|n|E|^2\|_{L^1}&\leq\|n\|_{L^2}\|E\|_{L^4}^2\\
&\leq\frac{1}{4}\|n\|_{L^2}^2+K^4(d)\|E_0\|_{L^2}^{4-d}\|\nabla E\|_{L^2}^d,
\end{split}\end{equation*}
\begin{equation*}\begin{split}\|B(E\otimes \overline{E})\|_{L^1}&\leq\eta\|E\otimes \overline{E}\|_{L^2}^2\\
&\leq\eta\|E\|_{L^4}^4\leq \eta K^4(d)\|E_0\|_{L^2}^{4-d}\|\nabla E\|_{L^2}^d.
\end{split}\end{equation*}
Consequently. we deduce from (5.4) that
\begin{equation}
\|\nabla E\|_{L^2}^2\leq H_2(0)+(1+\frac{\eta}{2})K^4(2)\|E_0\|_{L^2}^{2}\|\nabla E\|_{L^2}^2, \qquad\text{if}\; d=2,
\end{equation}
\begin{equation}
\|\nabla E\|_{L^2}^2\leq H_2(0)+(1+\frac{\eta}{2})K^4(3)\|E_0\|_{L^2}\|\nabla E\|_{L^2}^3, \qquad\text{if}\; d=3,
\end{equation}
By Lemma 5.2 to (5.5), (5.6) to yield (5.2) in the assumption (5.1) of Theorem 5.1.
\par
If k=2 and $d=2$, as the process of proof of Theorem 7.1 in \cite{La}, we have
\begin{equation*}
\frac{d}{dt}(m(t))\leq Cm(t)(1+\log m(t))
\end{equation*}
where $m(t)=\|\partial_tn\|_{L^2}^2+\|\nabla n\|_{L^2}^2+\|\partial_tE\|_{L^2}^2+1$, i.e.
\begin{equation}
\frac{d}{dt}(1+\log m(t))\leq C(1+\log m(t)).
\end{equation}
By the Gronwall lemma to (5.7) is given by
\begin{equation}
\|\partial_tn\|_{L^2}^2+\|\nabla n\|_{L^2}^2+\|\partial_tE\|_{L^2}^2+1\leq Ce^{e^{ct}}.
\end{equation}
Note that $\|\Delta E\|_{L^2}\leq C(\|E_t\|_{L^2}+\|\nabla n\|_{L^2}+1)$. Hence we obtain the result (5.3) as $k=2$.
\par
By mathematical induction, assume the
result (5.3) of Theorem 5.1 is valid for the case m = k +1. We now consider the case $m=k+2$,
applying the second equation of system (1.1) by the
operator $\partial^k$, taking the scalar product of $2\partial^kn$, integration by parts, we have
\begin{equation}\begin{split}
\frac{d}{dt}(\|\partial^kn_t\|_{L^2}^2&+\|\nabla\partial^kn\|_{L^2}^2)\leq2\|\partial^kn_t\|_{L^2}\|\Delta\partial^k|E|^2\|_{L^2}\\
&\leq2\|\partial^kn_t\|_{L^2}\left(\|E\|_{L^\infty}\|\Delta\partial^kE\|_{L^2}+\sum_{1\leq i,j\leq k+1}^{i+j=k+2}\|\partial^iE\partial^j\overline{E}\|_{L^2}\right)\\
&\leq Ce^{e^{\cdots^{e^{ct}}}}(\|\partial^kn_t\|_{L^2}^2+\|\Delta\partial^kE\|_{L^2}^2+1),
\end{split}\end{equation}
the last inequality is guaranteed by
\begin{equation*}\begin{split}
\|\partial^iE\partial^j\overline{E}\|_{L^2}^2&\leq\|\partial^{i-1}E\|_{L^2}^{1/2}\|\partial^{i+1}E\|_{L^2}^{1/2}
\|\partial^{j-1}E\|_{L^2}^{1/2}\|\partial^{j+1}E\|_{L^2}^{1/2}\\
&\leq C\|\Delta\partial^kE\|_{L^2}^{1/2}.
\end{split}\end{equation*}
Differentiating the first equation of system (1.1) with respect to the time variable, then applying the operator $\partial^k$,  Multiplying the resulting equation by $2\partial^k\overline{E}_t$, integration by parts, taking the imaginary part, it follows that
\begin{equation}
\frac{d}{dt}\|\partial^kE_t\|_{L^2}^2=Im2\int_{\mathbb{R}^2}\partial^k(nE)_t\partial^k\overline{E}_tdx-Im2\int_{\mathbb{R}^2}i\partial^k(E\otimes B)_t\partial^k\overline{E}_tdx.
\end{equation}
We first deal with the first term of right hand in (5.10) as follows
\begin{equation}\begin{split}
Im\int_{\mathbb{R}^2}\partial^k(nE)_t\partial^k\overline{E}_tdx=&Im\int_{\mathbb{R}^2}\partial^k(n_tE)\partial^k\overline{E}_tdx\\
&+Im\int_{\mathbb{R}^2}\sum_{0\leq j\leq k-1}^{i+j=k}\partial^in\partial^jE_t\partial^k\overline{E}_tdx,
\end{split}\end{equation}
 By virtue of the induction and interpolation inequality, we have
\begin{equation}\begin{split}
\|\partial^k(n_tE)\partial^k\overline{E}_t\|_{L^1}&\leq C\left(\|E\|_{L^\infty}\|\partial^kn_t\|_{L^2}+\sum_{0\leq j\leq k-1}^{i+j=k}\|\partial^in_t\partial^jE\|_{L^2}\right)
\|\partial^kE_t\|_{L^2}\\&
\leq C\left(\|\partial^kn_t\|_{L^2}+\sum_{0\leq j\leq k-1}^{i+j=k}\|\partial^in_t\|_{L^4}\|\partial^jE\|_{L^4}\right)\|\partial^kE_t\|_{L^2}\\&
\leq C\|\partial^kn_t\|_{L^2}\|\partial^kE_t\|_{L^2}+C\sum_{0\leq j\leq k-1}^{i+j=k}\|\partial^in_t\|_{L^2}^{1/2}\|\partial^{i+1}n_t\|_{L^2}^{1/2}\\&\qquad\times\|\partial^jE\|_{L^2}^{1/2}\|\partial^{j+1}E\|_{L^2}^{1/2}\|\partial^kE_t\|_{L^2}\\
&\leq Ce^{e^{\cdots^{e^{ct}}}}(\|\partial^kn_t\|_{L^2}^2+\|\partial^kE_t\|_{L^2}^2+1),
\end{split}\end{equation}
where $e^{e^{\cdots^{e^{ct}}}}$ denotes $(k+1)$-exponent, we have used the inequality $$\|E\|_{H^{k+1}}+\|n\|_{H^{k}}+\|\partial_tE\|_{H^{k-1}}+\|\partial_tn\|_{H^{k-1}}\leq Ce^{e^{\cdots^{e^{ct}}}},$$ which is guaranteed by induction assumption.
Similarly,
 \begin{equation}
\sum_{0\leq j\leq k-1}^{i+j=k}\|\partial^in\partial^jE_t\partial^k\overline{E}_t\|_{L^1}\leq C(\|\partial^{k+1}n\|_{L^2}^2+\|\partial^kE_t\|_{L^2}^2+1).
\end{equation}
Substituting (5.12) and (5.13) into (5.11) to yield
\begin{equation}
\left|Im\int_{\mathbb{R}^2}\partial^k(nE)_t\partial^k\overline{E}_tdx\right|\leq Ce^{e^{\cdots^{e^{ct}}}}(\|\partial^{k+1}n\|_{L^2}^2+\|\partial^kn_t\|_{L^2}^2+\|\partial^kE_t\|_{L^2}^2+1).
\end{equation}
Next, we will investigate the second term of right hand in (5.10). Note that $(iB)$
is real function, we have \begin{equation}\begin{split}
Im\int_{\mathbb{R}^2}i\partial^k(E\otimes B)_t\partial^k\overline{E}_t&dx=Im\int_{\mathbb{R}^2}i\partial^k(E\otimes B_t)\partial^k\overline{E}_tdx\\
&+Im\int_{\mathbb{R}^2}\sum_{0\leq l\leq k-1}^{l+j=k}(\partial^lE_t\otimes \partial^j(iB))\partial^k\overline{E}_tdx,
\end{split}\end{equation}
Since
\begin{equation}\begin{split}
|Im\int_{\mathbb{R}^2}&\sum_{0\leq l\leq k-1}^{l+j=k}(\partial^lE_t\otimes \partial^j(iB))\partial^k\overline{E}_tdx|\\&\leq
\sum_{0\leq l\leq k-1}^{l+j=k}\|\partial^lE_t\|_{L^4}\|\partial^jB\|_{L^4}\|\partial^k\overline{E}_t\|_{L^2}\\&\leq
\sum_{0\leq l\leq k-1}^{l+j=k}\|\partial^lE_t\|_{L^2}^{1/2}\|\partial^{l+1}E_t\|_{L^2}^{1/2}\|\partial^jB\|_{L^2}^{1/2}\|\partial^{j+1}B\|_{L^2}^{1/2}\|\partial^k\overline{E}_t\|_{L^2}
\\&\leq Ce^{e^{\cdots^{e^{ct}}}}(\|\partial^kE_t\|_{L^2}^{3/2}+\|\partial^kE_t\|_{L^2})\\&
\leq Ce^{e^{\cdots^{e^{ct}}}}(\|\partial^kE_t\|_{L^2}^{2}+1)
\end{split}\end{equation}
and
\begin{equation}\begin{split}
|Im\int_{\mathbb{R}^2}&i\partial^k(E\otimes B_t)\partial^k\overline{E}_tdx|\\&\leq
\left(\|E\|_{L^\infty}\|\partial^kB_t\|_{L^2}+\sum_{0\leq j\leq k-1}^{l+j=k}\|\partial^lE\|_{L^4}\|\partial^jB_t\|_{L^4}\right)\|\partial^k\overline{E}_t\|_{L^2}\\&\leq
\sum_{0\leq l\leq k-1}^{l+j=k}\|\partial^lE\|_{L^2}^{1/2}\|\partial^{l+1}E\|_{L^2}^{1/2}\|\partial^jB_t\|_{L^2}^{1/2}\|\partial^{j+1}B_t\|_{L^2}^{1/2}\|\partial^k\overline{E}_t\|_{L^2}\\&
\qquad\qquad+
 Ce^{e^{\cdots^{e^{ct}}}}\|\partial^k\overline{E}_t\|_{L^2}^2
\\&\leq Ce^{e^{\cdots^{e^{ct}}}}(\|\partial^kE_t\|_{L^2}^{3/2}+\|\partial^kE_t\|_{L^2}^2)\\&\leq Ce^{e^{\cdots^{e^{ct}}}}(\|\partial^kE_t\|_{L^2}^{2}+1),
\end{split}\end{equation}
where we have used the equality$$B=\frac{i\eta}{\Delta+\beta I}\nabla\times\nabla\times(E\otimes \overline{E}).$$
Inserting (5.16) and (5.17) into (5.15), it follows that
\begin{equation}
\left|Im\int_{\mathbb{R}^2}i\partial^k(E\otimes B)_t\partial^k\overline{E}_tdx\right|\leq  Ce^{e^{\cdots^{e^{ct}}}}(\|\partial^kE_t\|_{L^2}^{2}+1).
\end{equation}
Plugging (5.14) and (5.18) into (5.10) to deduce
\begin{equation}
\frac{d}{dt}\|\partial^kE_t\|_{L^2}^2\leq Ce^{e^{\cdots^{e^{ct}}}}(\|\partial^{k+1}n\|_{L^2}^{2}+\|\partial^kn_t\|_{L^2}^{2}+\|\partial^kE_t\|_{L^2}^{2}+1).
\end{equation}
Adding  (5.9) with (5.19), by the inequality
\begin{equation}\|\Delta\partial^kE\|_{L^2}\leq C(\|\partial^k\nabla n\|_{L^2}+\|\partial^kE_t\|_{L^2}+1),
\end{equation}
which is estimated by the first equation in system (1.1).
Hence we deduce
\begin{equation}\begin{split}
\frac{d}{dt}(\|\partial^kn_t\|_{L^2}^2&+\|\partial^k\nabla n\|_{L^2}^2+\|\partial^kE_t\|_{L^2}^2+1)\\&\leq Ce^{e^{\cdots^{e^{ct}}}}(\|\partial^kn_t\|_{L^2}^2+\|\partial^kE_t\|_{L^2}^2+\|\partial^k\nabla n\|_{L^2}^2+1).
\end{split}\end{equation}
By the Gronwall lemma to (5.21), using (5.20), we obtain
\begin{equation}
\|E\|_{H^{k+2}}^2+\|n\|_{H^{k+1}}^2+\|n_t\|_{H^{k}}^2+\|E_t\|_{H^{k}}^2\leq Ce^{e^{\cdots^{e^{ct}}}},
\end{equation}
where the $e^{e^{\cdots^{e^{ct}}}}$ denotes the $(k+2)$-exponent. The proof of Theorem 5.1 is completed.\hspace{\fill}$\blacksquare$
\begin{remark5} As $d=4$, if the initial data $(E,n,V)\in H^1\times L^2\times L^2$, $\|\nabla E_0\|_{L^2}^2\leq H_2(0)$ and $(1+\eta/2)CH_2(0)<1$, then system (1.1) in the case $(A_1)$ has a global weak solution
$$(E,n,V)\in \mathcal {C}(\mathbb{R}^+;H^1\times L^2\times H^{-1}),$$
where the constant $C$ satisfies the inequality
\begin{equation*}\|E\|_{L^4}\leq C\|\nabla E\|_{L^2}. \end{equation*}
The above proof is similar to the proof of the case $d=2,3$ in Theorem 5.1. In fact, we have
\begin{equation*}\begin{split}
\|\nabla E\|_{L^2}^2+&\frac{1}{4}\|n\|_{L^2}^2+\frac{1}{2}\|V\|_{L^2}^2\leq H_2(0)+(1+\frac{\eta}{2})\| E\|_{L^4}^4,\\&\leq
 H_2(0)+(1+\frac{\eta}{2})C\|\nabla E\|_{L^2}^4.
\end{split}\end{equation*}
Hence, it follows that
\begin{equation}
\|\nabla E\|_{L^2}^2\leq H_2(0)+(1+\frac{\eta}{2})C\|\nabla E\|_{L^2}^4.
\end{equation}
 By virtue of Lemma 5.2 to (5.23) yields the above result.
\end{remark5}
\begin{remark5}
There exists a family of self-similar blow-up  solution to the system (1.1) in 2D\cite{GGH, GM, GM1}. With assumption of small initial data in 2D,  we prove the global solution $(E,n,\partial_tn)\in \mathcal {C}(\mathbb{R}^+;H^k\times H^{k-1}\times H^{k-2})$ to system (1.1), $k\geq1$. In 3D, if the initial data is small enough, we obtain the unique and global solution $(E,n,\partial_tn)\in \mathcal {C}(\mathbb{R}^+;H^1\times L^2\times H^{-1})$. Moreover, for the 1D, system (1.1) becomes Eq.(1.2), the global well-posedness of solution $(E,n)\in L^2\times H^{-1/2}$, which obtianed by J. Colliander et. in \cite{CHT} is critial and optimal, because of the the ill-posedness of Eq.(1.2) in \cite{Ho}.
\end{remark5}
\textbf{Acknowledgments}

This work was partially supported by CPSF (Grant No.: 2013T60086) and NSFC (Grant No.: 11401122). The authors thank the references
for their valuable comments and constructive suggestions.\\


\begin{thebibliography}{99}
\small
\bibitem{Ad} H. Added and S. Added, Existence globale de
solutions fortes pour les \'{e}quations de la turbulence de Langmuir
en dimension 2, {\it C. R. Acad. Sci. Paris} {\bf 299} (1984)
551--554.
\bibitem{Ad1}
H. Added and S. Added, Equations of Langmuir turbulence and
nonlinear Schr\"{o}dinger equation: smoothness and approximation,
{\it J. Funct. Anal.} {\bf 79} (1988) 183--210.

\bibitem{B-C-D}
 H. Bahouri, J.Y. Chemin and R. Danchin, Fourier Analysis and Nonlinear Partial
Differential Equations. Springer--Verlag Berlin Heidelberg, 2011.

\bibitem{Bo}
J.M. Bony, Calcul symbolique et propagation des singularities\'{e}
pour les\'{e}quations aux d\'{e}riv\'{e}es partielles non
lin\'{e}aires. {\it Ann. Sci. \'{E}cole Norm. Sup.}  {\bf 14} (1981)
209--246.
\bibitem{B}
J. Bourgain and J. Colliander, On well-posedness of the Zakharov
system, {\it Int. Math. Res. Not.} {\bf 11} (1996) 515--546.

\bibitem{CCS}
J. Colliander, M. Czubak, and C. Sulem, Lower bound for the rate of blow-up of singular solutions of the Zakharov system in $\mathbb{R}^3$,
arXiv:1305.0324.
\bibitem{CHT}
L. Colliander, J. Holmer and N. Tzirakis, Low regularity global well-posedness for the
Zakharov and Klein-Gordon-Schr\"{o}dinger systems, {\it Trans. Amer. Math. Soc.} {\bf 360} (2008) 4619--4639.


\bibitem{GGH}
Zaihui Gan, Boling Guo and Daiwen Huang, Blow-up and nonlinear
instability for the magnetic Zakharov system, {\it J. Funct. Anal.}
{\bf 265} (2013) 953--982.

\bibitem{CW}
 T. Cazenave and F. Weissler, The Cauchy problem for the critical
nonlinear Schr\"{o}dinger equation in $H^S$, {\it Nonlinear Anal.
ATM} {\bf 14} (1990) 807--836.

\bibitem{GTV}
J. Ginibre, Y. Tsutsumi and G. Velo, On the Cauchy Problem for the
Zakharov System, {\it J. Funct. Anal.} {\bf 151} (1997) 384--436.

\bibitem{GM}
L. Glangetas and F. Merle, Existence of self-similar blow-up
solution for Zakharov equation in dimension two, Part I, {\it Comm.
Math. Phys.} {\bf 160} (1994) 173--215.

\bibitem{GM1}
 L. Glangetas and F. Merle,
Concentration properties of blow-up solutions and instability
results for Zakharov equation in dimension two, part II, {\it Comm.
Math. Phys.} {\bf 160} (1994) 349--389.

\bibitem{G}
R.T. Glassey, On the blowing-up of solutions to the Cauchy problem
for the nonlinear Schr\"{o}dinger equation, {\it J. Math. Phys.}
{\bf 18} (1977) 1794--1797.

\bibitem{Ho}
J. Holmer, Local well-posedness of the 1D Zakharov system, {\it Electronic J. Differ. Equations.} {\bf 24} (2007) 1--24.

\bibitem{KPV}
C. Kenig, G. Ponce and L. Vega, On the Zakharov and
Zakharov--Schulman systems, {\it J. Funct. Anal.} {\bf 127} (1995)
204--234.

\bibitem{KSH}
M.Kono, M.M. Skoric and D. Ter Haar, Spontaneous excitation of
magnetic fields and collapse dynamic ina Langmuir plasma, {\it J.
Plasma Phys.} {\bf 26} (1981) 123--146.

\bibitem{LPSS}
 M. Landman, G.C. Papanicolaou, C. Sulem and P.L. Sulem, Rate of
blow-up for solutions of the nonlinear Schr\"{o}dinger equation in
critical dimension, {\it Phys. Rev. A} {\bf 38} (1988) 3837--3843.

\bibitem{LPSSW}
 M. Landman, G.C. Papanicolaou, C. Sulem, P.L. Sulem and X.P. Wang,
 Stablility of isotropic self-similar dynamics for scalar collapse,
{\it Phys. Rev. A} {\bf 46} (1992) 7869--7876.

\bibitem{La}
C. Laurey, The Cauchy problem for a generalized Zakharov system,
{\it Differential Integral Equations} {\bf 8}  (1995) 105--130.

\bibitem{Me}
F. Merle, Blow-up results of virial type for Zakharov equations,
{\it Comm. Math. Phys.}  {\bf 175} (1996) 433--455.

\bibitem{Me1}
F. Merle, Lower bounds for the blowup rate of solutions of the
Zakharov equation in dimension two, {\it Comm. Pure Appl. Math.}
{\bf XLIX} (1996) 765--794.

\bibitem{Na}
H. Nawa, Asymptotic profiles of blow-up solutions of the nonlinear
Schr\"{o}dinger equation with critical power nonlinearity, {J. Math.
SOC. Japan} {\bf 46} (1994) 557--586.

\bibitem{OT}
T. Ozawa and Y. Tsutsumi, Existence and smoothing effect of
solutions for the Zakharov equations, {\it RIMS Kyoto Univ.} {\bf
28} (1992) 329--361.

\bibitem{PSSW}
G.C. Papanicolaou, C. Sulem, P.L. Sulem and  X.P. Wang, Singular
solutions of the Zakharov equations for Langmuir turbulence, {\it
Phys. Fluids B} {\bf 3} (1991) 969--980.

\bibitem{Su}
C. Sulem and P.L. Sulem, Quelques r\'{e}sultats de
r\'{e}gularit\'{e} pour les \'{e}quations de la turbulence de
Langmuir, {\it C. R. Acad. Sci. Paris} {\bf 289} (1979) 173--176.

\bibitem{We}
M.I. Weinstein, Nonlinear Schr\"{o}dinger Equations
and Sharp Interpolation Estimates, {\it Comm. Math. Phys.} {\bf 87} (1983) 567--576.
\bibitem{Za}
V.E. Zakharov, Collapse of Langmuir waves, {\it Soviet Phys.} JETP,
{\bf 35} (1972) 908--914.
\end{thebibliography}
\end{document}